\documentclass[11pt]{amsproc}

\usepackage{amssymb,a4}
\renewcommand{\labelenumi}{(\roman{enumi})}

\newcommand{\bbfC}{{\Bbb C}}

\newcommand{\bbfP}{{\Bbb P}}
\newcommand{\bbfR}{{\Bbb R}}
\newcommand{\bbfZ}{{\Bbb Z}}

\newcommand{\PP}{\bbfC\bbfP^2}
\newcommand{\PQ}{\bbfP^1\times\bbfP^1}
\newcommand{\hol}{{\cal O}}
\newcommand{\im}{\mbox{\rm im\,}}                   
\newcommand{\oo}{\mbox{\rm o}}
\newcommand{\Pic}{\mbox{\rm Pic\,}}                 

\newcommand{\fun}[1]{\textstyle -\frac{#1}{2}K}     
\newcommand{\fund}{\fun{1}}                         
\newcommand{\fb}[1]{\textstyle K^{-\frac{#1}{2}}}   
\newcommand{\fdb}{\fb{1}}                           
\newcommand{\fdbd}{\textstyle K^{\frac{1}{2}}}      
\newcommand{\can}[1]{\textstyle K^{-#1}_S}          
\newcommand{\canS}{\can{1}}                         
\newtheorem{thm}{Theorem}[section]
\newtheorem{lemma}[thm]{Lemma}
\newtheorem{prop}[thm]{Proposition}
\newtheorem{cor}[thm]{Corollary}
\theoremstyle{definition}
\newtheorem{rem}[thm]{Remark}
\newtheorem{defi}[thm]{Definition}
\theoremstyle{remark}

\newenvironment{main}[1]{{\par\bf #1} \sl}{\par}

\newcommand{\pf}{{\sc Proof:\quad}}
\renewcommand{\bf}{\bfseries}
\renewcommand{\sl}{\slshape}
\newcommand{\cal}{\mathcal}
\renewcommand{\Bbb}{\mathbb}

\begin{document}

\title[Algebraic dimension of twistor spaces]{On the algebraic dimension of
  twistor spaces over the connected sum of four complex projective planes} 
\author { B.\ Kreu{\ss}ler}
\address{  FB Mathematik, Universit\"at Kaiserslautern\\
  D--67\,653 Kaiserslautern, Germany }
\email{kreusler@mathematik.uni-kl.de}

\date{July 1996}


\begin{abstract}
\noindent
We study the algebraic dimension of twistor spaces of positive type over
$4\PP$. We show that such a twistor space is Moishezon if and only if its
anticanonical class is not nef. More precisely, we show the equivalence of
being Moishezon with the existence of a smooth rational curve having negative
intersection number with the anticanonical class. Furthermore, we give precise
information on the dimension 
and base locus of the fundamental linear system $|\fund|$.  This implies, for
example, $\dim|\fund|\leq a(Z)$. We characterize those twistor spaces over
$4\PP$, which contain a pencil of divisors of degree one by the property
$\dim|\fund| = 3$.  
\end{abstract}

\maketitle
\section{Introduction}
\label{intro}

Twistor spaces usually arise in four--dimensional conformal geometry. Their
construction reflects the impossibility to equip in general a four--dimensional
conformal manifold $M$ with a compatible complex structure. It was shown in
\cite{AHS} that the conformal metric on $M$ is self--dual if and only if the
twistor space $Z$ associated to $M$ carries, in a natural way, the structure of
a complex manifold. Therefore, the conformal geometry of $M$ is closely related
to the holomorphic geometry of $Z$. Since we shall only work with methods of
complex geometry, we can use the following definition:
\\ A twistor space $Z$ is a complex three--manifold with the following
additional structure: 
\begin{itemize}
  \item
  a proper differentiable submersion $\pi:Z\rightarrow M$ onto a real
  differentiable four--manifold $M$. The fibres of $\pi$ are holomorphic curves
  in $Z$ being isomorphic to $\bbfC\bbfP^1$ and having normal bundle in $Z$
  isomorphic to $\hol(1)\oplus\hol(1)$;
  \item
  an anti--holomorphic fixed point free involution $\sigma:Z\rightarrow Z$
  with $\pi\sigma=\pi$.
\end{itemize}
The fibres of $\pi$ are called ``real twistor lines'' and the involution
$\sigma$ is called the ``real structure''. A geometric object will be called
``real'' if it is $\sigma$--invariant. For example, a line bundle ${\cal L}$ on
$Z$ is real if $\sigma^\ast\bar{{\cal L}}\cong {\cal L}$, and a complex
subvariety $D\subset Z$ is real if $\sigma(D) = D$. Instead of $\sigma(D)$ we
shall often write $\bar{D}$.

We only consider compact and simply connected twistor spaces.

At the beginning of the 80's, the first classification result emerged in
\cite{FK}, \cite{Hit2}:\\
There exist exactly two compact K{\"a}hlerian twistor spaces. They are
automatically projective algebraic. The corresponding Riemannian
four--manifolds 
are the $4$--sphere $S^4$ and the complex projective plane $\PP$ (with
Fubini--Study metric).
This was generalized in \cite{Camsc} to the result that a twistor space which
is bimeromorphic to a compact K{\"a}hler manifold must be Moishezon and simply
connected. This implies (see \cite{Don}, \cite{F}) that $M$ is homeomorphic to
the connected sum $n\PP$ for some $n\geq 0$.

New examples of Moishezon twistor spaces were constructed by Y.S.\ Poon
\cite{Po1} (case $n = 2$) and C.\ LeBrun \cite{LeB}, H.\ Kurke \cite{Ku} (case
$n\geq 3$). 

Nowadays the situation is well understood for $n\leq 3$ (cf. \cite{Hit2},
\cite{FK}, \cite{Po1}, \cite{KK}, \cite{Po2}). To become more precise, we have
to introduce the notion of the ``type'' of a twistor space. By a result of
R.~Schoen \cite{Sch}, every conformal class of a compact Riemannian
four--manifold contains a metric of constant scalar curvature. Its sign will be
called the {\sl type\ } of the twistor space. This is an invariant of the
conformal class, hence of the twistor space.

It was shown in \cite{Po3} that a Moishezon twistor space is always of positive
type. If $n \leq 3$, the converse is also true.

In this paper we focus on the positive type case, for two reasons. One reason
is that we can then apply Hitchin's vanishing theorem (\ref{Hit}). The 
other reason is the following: a result of P.\ Gauduchon \cite{Gau} implies
that any twistor space of negative type has algebraic dimension zero. From the
results of M.\ Pontecorvo \cite{Pon} we easily derive that a twistor space of
type zero over $n\PP$ must also have algebraic dimension zero. It is not clear
whether there exist twistor spaces of non--positive type over $n\PP$.

Computation of algebraic dimension is, therefore, interesting only in the case
of positive type. A very important tool to compute the algebraic dimension of
twistor spaces is the result of Y.S.\ Poon \cite{Po3} (see also \cite{Pon})
stating that the algebraic dimension is equal to the Iitaka dimension of the
anticanonical bundle (cf. Section \ref{prelim}).

From \cite{DonF} and \cite{Cam}, \cite{LeBP} it is known that the generic
twistor space over $n\PP$ has algebraic dimension one (if $n = 4$),
respectively zero (if $n\geq 5$).

For the case $n=4$, the characterizing property $c_1^3 = 0$ is of central
importance.

In this paper we study the following
\begin{main}{Problem:}
Compute the algebraic dimension $a(Z)$ of a twistor space $Z$ over $4\PP$ in
terms of geometric or numeric properties of certain divisors on $Z$.
\end{main}
A first attempt to tackle this problem was made by Y.S.\ Poon \cite[Section
7]{Po2}. He assumes, additionally, the existence of a divisor $D$ of degree one
on $Z$. He studies a birational map $D\rightarrow\bbfP^2$, which is the
blow--up of four points. He seems to assume that these four points are
actually in $\bbfP^2$ (no infinitesimally near blown--up points). If these four
points are in a special position he obtains $a(Z) = 3$. In the case of general
position he can only show: $a(Z) \leq 2$.

We shall, in general, not assume the existence of divisors of degree
one. Because in case $n=4$ there exists at least a pencil of so--called
fundamental divisors, 
we shall study their geometry to obtain our results. If $S\subset Z$ is a real
fundamental divisor, we have a birational map $S\rightarrow\PQ$ which is the
blow--up of eight points. We shall study in detail the possible positions for
these points. We take into account that some of these points can be
infinitesimally near each other. We are able to derive the algebraic dimension
$a(Z)$ from the knowledge of the positions of these eight points.

Similar considerations were made for general $n\geq 4$ in the paper \cite{PP}.
But the authors of that paper are intersted in  a study of small
deformations of well--known Moishezon twistor spaces, and so they investigate
only the case without infinitely near blown--up points.

As a consequence of our results, we give a new characterization of the
twistor spaces over $4\PP$ which are first described by C.\ LeBrun \cite{LeB}
(with methods from differential geometry). From the point of 
view of complex geometry the twistor space structure on these complex manifolds
was found by H.\ Kurke \cite{Ku}. Following the literature, we call them {\em
  LeBrun twistor spaces}.
These twistor spaces are characterized in \cite{Ku} and \cite{Po2} by the
property to contain a pencil of divisors of degree one. In the case 
$n=4$ we show (Theorem (\ref{cb})) that they can also be characterized by the
property $h^0(\fdb) = 4$ or by the structure of the base locus of $|\fund|$.

Besides this, our main results are a precise description of the set of
irreducible curves intersecting $\fdb$ negatively (Theorem \ref{ncurves}) and 
the following theorems, where $Z$ denotes
always a simply connected compact twistor space of positive type over $4\PP$:

\setcounter{section}{6}
\setcounter{thm}{1}
\begin{thm}
  $a(Z) = 3 \iff \fdb$ is not nef;\\[-0.8mm]
  $a(Z) = 2 \iff \fdb$ is nef and $\exists m\geq 1: h^1(\fb{m}) \ne 0$;%
\\[-0.8mm]
  $a(Z) = 1 \iff \forall m\geq 1: h^1(\fb{m}) = 0$.
\end{thm}

\begin{thm}
  The following conditions are equivalent:
  \begin{enumerate}
  \item $a(Z) = 3$;\vspace*{-1mm}
  \item $\fdb$ is not nef;\vspace*{-1mm}
  \item there exists a smooth rational curve $C\subset Z$ with $C.(\fund) < 0$.
  \end{enumerate}
\end{thm}

\setcounter{thm}{5}
\begin{thm}
  $a(Z) \geq \dim|\fund|$.
\end{thm}

\begin{thm}
  If $\dim|\fund|\geq 2$, then:\\
  $a(Z) = 2 \iff \fdb$ is nef $\iff |\fund|$ does not have base points.
\end{thm}

\setcounter{section}{1}
\setcounter{thm}{0}
This paper is organized as follows:\\
In Section \ref{prelim} well--known but necessary facts about simply connected
compact twistor spaces of positive type are collected.

Also Section \ref{eins} has preparatory character. We study there the structure
of fundamental divisors for general $n$, using results and techniques contained
in \cite{PP}. Technically important for the following sections will be
Proposition \ref{types} where the structure of effective anticanonical curves
on real fundamental divisors is described in detail.

In the remaining three sections we assume $n=4$.

In Section \ref{zwei} we study the case where the anti--canonical bundle
$K_Z^{-1}$ is nef (in the sense of Mori theory). We shall prove that the
algebraic dimension is, in this case, at most two. We also see how to
distinguish between algebraic dimension one and two. This generalizes results
of \cite{CK}. 

In Section \ref{drei} we assume $K_Z^{-1}$ to be not nef. We 
collect detailed information on the fundamental linear
system $|\fund|$ and on the set of curves which intersect $\fdb$ negatively. In
this cases the algebraic dimension is three.  

The final Section \ref{vier} combines the results of the previous part to
prove the main theorems stated above.

\begin{main}{Acknowledgement}
  I thank Fr\'ed\'eric Campana for encouragement and stimulating discussions.
\end{main}
%
\section{Preliminaries}
\label{prelim}

We briefly collect well--known facts which will be frequently used later. We
refer the reader to \cite{AHS}, \cite{ES}, \cite{Hit2}, \cite{Kr}, \cite{Ku}
and \cite{Po1}. For brevity, we assume, throughout this section, $Z$ to be a
simply connected compact twistor space of positive type. As mentioned in
Section \ref{intro} the corresponding Riemannian four--manifold $M$ is
homeomorphic 
to $n\PP$.

\subsubsection*{Cohomology ring of $Z$}

$H^i(Z,\bbfZ)$ is a free $\bbfZ$--module.\\
$H^1(Z,\bbfZ) = H^3(Z,\bbfZ) = H^5(Z,\bbfZ) = 0$ and $H^0(Z,\bbfZ) \cong
H^6(Z,\bbfZ) \cong \bbfZ$. \\
$H^2(Z,\bbfZ)$ and  $H^4(Z,\bbfZ)$ are free modules of rank $n+1$. There exists
a basis $x_1,\dots, x_n, w$ of $H^2(Z,\bbfZ)$ such that the pull--back
$H^2(Z,\bbfZ) \rightarrow H^2(F,\bbfZ)\cong \bbfZ$ (for any real twistor line
$F \subset Z$) sends $x_i$ to $0$ and $w$ to the positive generator.

The cohomology ring $H^\ast(Z,\bbfZ)$ is isomorphic to the graded ring $\bbfZ
[x_1,\dots, x_n, w]/R$ where $R$ is the ideal generated by 
\[ x_i^2 - x_j^2,\quad x_ix_j\; (i\ne j),\quad w^2 + w\sum_{i=1}^n x_i +
x_1^2.\] 
The grading is given by $\deg x_i = \deg w = 2$.\\
$H^4(Z,\bbfZ)$ is a free $\bbfZ$--module with generators $wx_1,\dots, wx_n,
w^2$. The dual class of  a real twistor fibre $F\subset Z$ is $-x_i^2 \in
H^4(Z,\bbfZ)$. 

$c_1(Z) = 4w + 2\sum_{i=1}^n x_i, \quad c_2(Z) = -6x_1^2 = 6F, \quad c_3(Z) =
2(n + 2)$. This yields the following Chern numbers:
$c_1^3 = 16(4 - n), \quad c_1c_2 = 24, \quad c_3 = 2(n + 2)$.

\subsubsection*{Cohomology of sheaves}

The main reason to assume $Z$ to be of positive type is Hitchin's vanishing
theorem. We shall only use the following special case:
\begin{thm}[Hitchin \cite{Hit}]\label{Hit}
  If $Z$ is of positive type then we have for any ${\cal L}\in \Pic(Z)$
  \begin{eqnarray*}
  \deg({\cal L})\le -2 &\Rightarrow& H^1(Z,{\cal L})=0.
  \end{eqnarray*}
\end{thm}

On the other hand, since the twistor lines cover $Z$, we obtain:
\begin{eqnarray*}
\deg({\cal L})\le -1 &\Rightarrow& H^0(Z,{\cal L})=0.
\end{eqnarray*}
By Serre duality this gives the following important vanishing results:
\begin{eqnarray}
\deg({\cal L})\ge -2 &\Rightarrow& H^2(Z,{\cal L})=0,\\
\deg({\cal L})\ge -3 &\Rightarrow& H^3(Z,{\cal L})=0.
\end{eqnarray}

In particular, we obtain $h^2(\hol_Z) = h^3(\hol_Z) = 0$. Because $Z$ is simply
connected, we also have $h^1(\hol_Z) = 0$. Hence, we obtain an isomorphism of
abelian groups, given by the first Chern class:
\[\Pic(Z)\stackrel{\sim}{\rightarrow}H^2(Z,\bbfZ).\]
There exists a unique line bundle whose first Chern class is $\frac{1}{2}c_1$.
We shall denote it by $\fdb$. Following Poon, we call it the {\sl fundamental}
line bundle.  The divisors in the linear system $|\fund|$ will be called {\sl
  fundamental divisors}. The description of the cohomology ring gives
$(\fund)^3 = 2(4 - n)$. If $S\in|\fund|$ is a smooth fundamental divisor, we
obtain by the adjunction formula $\canS\cong\fdb\otimes\hol_S$. If $n\leq 4$,
there exist smooth real fundamental divisors (cf.\ \cite[Lemma 3.1]{CK}).

The degree of a line bundle ${\cal L}\in \Pic(Z)$ will be by definition the
degree of its restriction to a real twistor line. For example,
$\deg(\fdb)=2$. We obtain in this way a {\sl surjective\/} degree map
\[\deg :\Pic(Z)\twoheadrightarrow \bbfZ.\]
From the above equations on Chern numbers we obtain, by applying the
Riemann--Roch theorem,
\begin{equation}\label{RR}
\chi(Z,\fb{m})=m+1+2(4-n){\binom{m+2}{3}}.
\end{equation}

\subsubsection*{Algebraic dimension}

We denote by $a(Z)$ the algebraic dimension of $Z$, which is by definition the
transcendence degree of the field of meromorphic functions of $Z$ over $\bbfC$.
If $\dim Z = a(Z)$, then $Z$ is called Moishezon.  To compute the algebraic
dimension of twistor spaces we shall frequently use, without further reference,
the following theorem of Y.S.\ Poon:
\begin{thm}{\cite{Po3}, \cite[Prop.\ 3.1]{Pon}}\\
 $ \begin{array}[t]{lcl}
    \kappa(Z, K^{-1}) \geq 0& \Rightarrow& a(Z) = \kappa(Z, K^{-1})\\
    \kappa(Z, K^{-1}) = -\infty& \Rightarrow& a(Z) = 0.
  \end{array}$  
\end{thm}

The number $\kappa(Z, K^{-1})$ is usually called the Iitaka dimension (or
L--dimension = line bundle dimension) of the line bundle $K^{-1}$. Its
definition generalizes the well--known notion of Kodaira dimension. For
details, including the following facts, we refer the reader to \cite{U}. \\ For
any line bundle ${\cal L}\in\Pic(Z)$ there holds: $\dim Z\geq a(Z)\geq
\kappa(Z, {\cal   L}) $.\\ 
If $f:Z\rightarrow Y$ is a dominant morphism, then $a(Z) \geq
a(Y)$.  Particularly, if $f:Z\rightarrow \bbfP^N$ is a meromorphic map, then
$a(Z) \geq \dim f(Z)$, because any projective variety is Moishezon.\\ 
If we define 
$g:=\gcd\{m\in \bbfZ\; |\; m > 0, h^0(Z, L^m) \ne 0\}$ and denote by
$\Phi_{|L^m|}$ the meromorphic map given by the linear system $|L^m|$, then
$\kappa(Z, L) = \max\{\dim \Phi_{|L^m|}(Z)\;|\;m \in g\bbfZ, \; m>0\}$. If
there exists a polynomial $P(X)$ such that for all large positive $m\in \bbfZ$
we have $h^0(Z, L^{mg}) \leq P(m)$, then $\kappa(Z, L)\leq \deg P$.

We apply these basic facts to obtain our first result on the algebraic
dimension in case $n = 4$. The following proposition is a generalization of a
result contained in \cite{CK}. For convenience we introduce the following
\begin{main}{Definition:}
  If there exists an integer $m\ge 1$ with $h^1(\fb{m})\ne 0$ then we define
  $\tau := \min\{m | m\ge 1, h^1(\fb{m})\ne 0\}$. Otherwise we set $\tau :=
  \infty$.
\end{main}

\begin{prop}
\label{tau}
Let $Z$ be a simply connected compact twistor space of positive type with
$c_1^3 = 0$. Then:\\
\begin{tabular}[t]{rl}
  (i) & $a(Z)\ge 1$\\
  (ii) & $a(Z) = 1 \iff \forall m\ge 1\quad h^1(\fb{m})=0$.
\end{tabular}
\end{prop}

\pf 
From Riemann--Roch and Hitchin's vanishing theorem we know:
$h^0(\fb{m})=m+1+h^1(\fb{m})$. Therefore, $a(Z)=\kappa(Z,\fdb)\ge 1$ and if
$\tau=\infty$ we have $\kappa(Z,\fdb)=1$.

Assume $\tau<\infty$ and $a(Z)=1$. Let $S\in|\fund|$ be smooth and real. (Such
a divisor exists, because we assume $n=4$, cf.\ \cite{CK}.) Since
$h^1(\fb{\tau-1})=0$ the exact sequence
\[0 \rightarrow\fb{\tau-1} \rightarrow \fb{\tau} \rightarrow \can{\tau} 
\rightarrow 0\]
gives an exact sequence
\[
0 \rightarrow H^0(\fb{\tau-1}) \rightarrow H^0(\fb{\tau}) \rightarrow
H^0(\can{\tau}) \rightarrow 0.
\]
Since $h^1(\fb{\tau})\ge1$ we have, furthermore,
$h^0(\fb{\tau})=\tau+1+h^1(\fb{\tau})\ge\tau+2$. The linear system
$|\fun{\tau}|$ cannot have a fixed component since $\tau S\in |\fun{\tau}|$ and
$\dim |S| = \dim |\fund| \ge 1$. If necessary blow up $Z$ to obtain a morphism
$\Phi_\tau:\tilde{Z}\rightarrow \bbfP^d$ defined by $|\fun{\tau}|$. Here
$d:=\dim |\fun{\tau}| \ge \tau+1\ge 2$. By assumption $\dim
\Phi_\tau(\tilde{Z})=1$. Since the curve $\Phi_\tau(\tilde{Z})$ is not
contained in a linear subspace of $\bbfP^d$, its degree must be at least $d$.
Hence, a generic member of the linear system $|\fun{\tau}|$ is the sum of
$\lambda$ algebraically equivalent divisors and so it is linearly equivalent to
$\lambda S_0$ with $\lambda\ge d\ge\tau+1$. This gives
$2\tau=\deg(\fun{\tau})=\lambda\deg(S_0),$ which is only possible if $\lambda =
2\tau$ and $\deg(S_0)=1$. But then we have infinitely many divisors of degree
one in $Z$. This implies $a(Z)=3$ by the Theorem of Kurke--Poon (see \cite{Ku},
\cite{Po2}). This contradiction proves the proposition.\qed

\begin{rem}
  If $|-K_S|$ contains a smooth curve $C$, then we computed in \cite{CK} that
  $\tau$ is the order of $N:=\canS\otimes\hol_C$ in the Picard group $\Pic{C}$
  of the elliptic curve $C$. Under this additional assumption Proposition
  \ref{tau} was shown in \cite{CK}.
\end{rem}

%
\section{The structure of fundamental divisors}
\label{eins}

In this section $Z$ always denotes a simply connected compact twistor space of
positive type.
\begin{lemma}\label{pic}
  Let $S\in|\fund|$ be a smooth surface. Then the restriction map $\Pic{Z}
  \rightarrow \Pic{S}$ is injective.

\end{lemma}
\pf
By assumption we have $h^1(\hol_Z)=h^2(\hol_Z)=0$. Since $S$ is a rational
surface \cite{Po1}, we also have $h^1(\hol_S)=h^2(\hol_S)=0$. Therefore, taking
the first Chern class defines isomorphisms $\Pic{Z}
\stackrel{\sim}{\rightarrow} H^2(Z,\bbfZ)$
and $\Pic{S} \stackrel{\sim}{\rightarrow} H^2(S,\bbfZ)$. Let us denote the
inclusion of $S$ into $Z$ by $i$. The above isomorphisms transform then the
restriction morphism $\Pic{Z} \rightarrow \Pic{S}$ into the map $i^\ast$ on
cohomology groups. 

We shall apply standard facts from algebraic topology to verify the injectivity
of $i^\ast$. Let $\oo_S \in H_4(Z,\bbfZ)$ and $\oo_Z \in H_6(Z,\bbfZ)$ be the
fundamental classes of $S$ and $Z$ respectively. By $d_Z(S) \in H^2(Z,\bbfZ)$
we denote the Poincar\'e dual of $i_\ast(\oo_S) \in H_4(Z,\bbfZ)$, this means
$i_\ast(\oo_S) = d_Z(S) \smallfrown \oo_Z$ (cap--product).

For any cohomology class $\alpha \in H^2(Z,\bbfZ)$ we obtain by the
associativity of cap--product $\alpha \smallfrown i_\ast(\oo_S) = \alpha
\smallfrown (d_Z(S) \smallfrown \oo_Z) = (\alpha \smallsmile d_Z(S))
\smallfrown \oo_Z$. The naturalness of cap--product implies $\alpha \smallfrown
i_\ast(\oo_S) = i_\ast(i^\ast(\alpha) \smallfrown \oo_S)$. Therefore, we obtain
a commutative diagram:

\centerline{
$\begin{array}[t]{ccc} 
  H^2(Z,\bbfZ) & \stackrel{{\scriptstyle i^\ast}}{\longrightarrow} &
  H^2(S,\bbfZ)\\ 
  & & \downarrow {\scriptstyle\smallfrown \oo_S}\\
  \Bigg\downarrow {\scriptstyle \smallsmile d_Z(S)} & & H_2(S,\bbfZ)\\
  & & \downarrow {\scriptstyle i_\ast}\\
  H^4(Z,\bbfZ) & \stackrel{{\scriptstyle \smallfrown \oo_Z}}{\longrightarrow} &
  H_2(Z,\bbfZ)
\end{array}$}

Since, by Poincar\'e duality, the cap--product with $\oo_Z$ is an isomorphism,
we obtain $\ker(i^\ast) \subset \ker( \smallsmile d_Z(S))$.

The description of the cohomology ring given above allows us to compute the
kernel of the cup--product with the dual class $d_Z(S)$ of S.  With the
notation of Section \ref{prelim} the elements $x_1,\dots, x_n,\omega$ form a
basis of the free $\bbfZ$--module $H^2(Z,\bbfZ)$. The dual class of $S$ is
$d_Z(S) = c_1(\fdb) = 2\omega + x_1 + \dots + x_n$. If we use $\omega
x_1,\dots,\omega x_n,x_1^2$ as basis of $H^4(Z,\bbfZ)$ then the cup--product
with $d_Z(S)$ is described by the $(n+1)\times(n+1)$--matrix:
\[ \left(
\begin{array}{rrrrr}
  2&0&\dots&0&1\\ 0&2&&0&1\\ \vdots& & &\vdots &\vdots\\ 0& &2 &0&1\\ 
  0&\dots&0&2&1\\ -1&\dots&-1&-1&-2
\end{array}
\right)\] 
whose determinant is equal to $2^{n-1}(n-4)$. 

If $n\ne 4$ we obtain the injectivity of the map $\alpha \mapsto \alpha
\smallsmile d_Z(S)$ and thus of the restriction map $\Pic{Z} \hookrightarrow
\Pic{S}$. If $n = 4$ it is easy to see that $\alpha \smallsmile d_Z(S) = 0$ if
and only if $\alpha \in \bbfZ\cdot d_Z(S) \subset H^2(Z,\bbfZ)$. To prove the
injectivity of $\Pic{Z} \rightarrow \Pic{S}$ it remains, therefore, to show
that $\fb{m}\otimes\hol_S \cong \hol_S$ implies $m = 0$.

By adjunction we have $\fb{m}\otimes\hol_S \cong \can{m}$. But $S$ is rational,
hence $\Pic{S}$ is torsion free and $K_S\ncong \hol_S$. Thus,
$\can{m}\cong\hol_S$ if and only if $m=0$. This proves the Lemma. \qed
 
\begin{lemma}\label{conn}
  Let $Z$ be a simply connected compact twistor space of positive type and
  $D\subset Z$ a divisor of degree one. If $S\in|\fund|$ is a smooth surface,
  then $C:=D\cap S$ is connected.
\end{lemma}
\pf
We shall show $h^0(\hol_C)=1$, which implies connectedness of $C$.

Consider first the exact sequence
\begin{equation}
  \label{seqelm}
  0 \rightarrow \hol_Z(-\bar{D}) \rightarrow \hol_Z \rightarrow \hol_{\bar{D}}
  \rightarrow 0.
\end{equation}
From $h^1(\hol_Z)=0$ we obtain
$h^1(\hol(-\bar{D}))=h^0(\hol_{\bar{D}}) - h^0(\hol_Z)=0$ since $Z$ and
$\bar{D}$ are connected.

As $D+\bar{D}\in |\fund|$ we obtain an exact sequence 
\[ 0 \rightarrow K(\bar{D}) \rightarrow \fdbd  \rightarrow \hol_D(-C)
\rightarrow 0. \] But the degree of $\fdbd$ is $-2$ and, therefore, Hitchin's
vanishing theorem gives $h^i(\fdbd)=0$ for all $i$. Therefore, using 
Serre duality, $h^1(\hol_D(-C))= h^2(K(\bar{D})) = h^1(\hol_Z(-\bar{D})) = 0$.

Consider finally the exact sequence
\[ 0 \rightarrow \hol_D(-C) \rightarrow \hol_D  \rightarrow \hol_C \rightarrow
0. \] We have $h^0(\hol_D)=1$ since any divisor of degree one is connected
\cite{Po1}. Because $C$ is effective, $h^0(\hol_D(-C))$ must vanish. Hence,
$h^0(\hol_C)=h^0(\hol_D)=1$. \qed

\begin{lemma}[cf.\ \cite{PP}, p.\ 693]\label{noreal}
  Let $Z$ be as above and $S\in|\fund|$ an irreducible real divisor. Then $S$
  is smooth and contains a real twistor fibre $F\subset S$. The linear system
  $|F|$ is one--dimensional and its real elements are precisely the real
  twistor fibres contained in $S$.
\end{lemma}

\pf
The smoothness of $S$ was shown in \cite[Lemma 2.1]{PP1}. If $S$ does not
contain a real twistor fibre, the restriction of the twistor fibration to $S$
would give an unramified double cover over a simply connected manifold, since
$Z$ does not contain real points. But $S$ is connected and must, therefore, 
contain a real twistor fibre $F$. From the adjunction formula we obtain $F^2 =
0$ on $S$. Hence, we have an exact sequence 
$0 \rightarrow \hol_S \rightarrow \hol_S(F) \rightarrow \hol_F \rightarrow 0$.
From $h^1(\hol_S) = 0$ we infer, therefore, $\dim |F| = 1$.

Since the linear system $|F|$ defines a flat family of curves in $S$, its
elements form a curve in the Douady space ${\cal D}$ of curves on $Z$ (cf.\ 
\cite{Dou}). Since $h^0(\bbfP^1,\hol(1)\oplus\hol(1))=4$ and
$h^1(\bbfP^1,\hol(1)\oplus\hol(1))=0,\;\;{\cal D}$ is a four--dimensional
complex manifold near points which correspond to smooth rational curves on $Z$
with normal bundle $\hol(1)\oplus\hol(1)$. The real structure of $Z$ induces
one on ${\cal D}$.  If the set of real points ${\cal D}(\bbfR)$ is non--empty,
then it is a four--dimensional real manifold near points as before.  Since the
real twistor lines are smooth rational curves with the above normal bundle, the
real manifold $M=4\PP$ is a submanifold of ${\cal D}(\bbfR)$. Since $M$ is
compact and has the same dimension as ${\cal D}(\bbfR)$, it must be a connected
component.

The set $U$ of members of $|F|$ which are smooth rational curves with normal
bundle $\hol(1)\oplus\hol(1)$ is open and dense in $\bbfP^1\cong|F|$ with
respect to the Zariski topology. Therefore, the set $U(\bbfR)$ of real points
in $U$ is open and dense in the one--sphere of real members of $|F|$. Since
${\cal D}(\bbfR)$ is smooth near $M$ and $M$ is a component of ${\cal
  D}(\bbfR)$, we have $U(\bbfR)\subset M$. But $M$ is compact and must,
therefore, contain the closure of $U(\bbfR)$ in ${\cal D}(\bbfR)$, which is the
set of all real members of $|F|$. Therefore, any real member of $|F|$ is a real
twistor fibre and, in particular, smooth and ireducible. This proves the
claim.\qed

To obtain more information on the structure of real irreducible fundamental
divisors $S\in|\fund|$ one can study the morphism $S\rightarrow \bbfP^1$ given
by $|F|$ (cf.\ \cite[p.\ 693]{PP}).
Since the general fibre of this morphism is a smooth rational curve it factors
through a rational ruled surface. Since $(-K_S)^2 = (\fund)^3 = 8-2n$, the
surface $S$ is a blow--up of  a ruled surface at $2n$ points. The exceptional
curves of these blow--ups are contained in fibres of the morphism $S\rightarrow
\bbfP^1$. By Lemma \ref{noreal} none of the exceptional curves is real and none
of the blown--up points lie on a real fibre of the ruled surface. Using this,
in \cite[Lemma 3.5]{PP} it has been shown that the ruled surface is isomorphic
to $\PQ$. Therefore, we obtain a morphism $\sigma:S\rightarrow\PQ$ which is a 
succession of blow--ups. Let us equip $\PQ$ with the real structure given by
the antipodal map  on the first factor and the  usual real structure on the
second. Then $\sigma$ is equivariant (or ``real''). Since we can always
contract a conjugate pair of disjoint $(-1)$--curves, $\sigma$ is the
succession of $n$ blow--ups. At each step a conjugate pair of points is
blown--up to give a surface without real points.

We should bear in mind that it is possible to have
infinitesimally near blown--up points.  

As in \cite{PP} we shall call curves of type $(1,0)$ on $\PQ$ ``lines'' and
curves of type $(0,1)$ ``fibres''. Then there do not exist real lines. But
the images of real twistor fibres in $|F|$ are exactly the real ``fibres''.

\begin{lemma}\label{smocomp}
  Equip $\PQ$ with the real structure $((a_0:a_1),(b_0:b_1))\mapsto
  ((\bar{a}_1:-\bar{a}_0),(\bar{b}_0:\bar{b}_1))$ as described above. Then the
  reduced components of any \underline{real} member of $|\hol(2,2)| =
  |-K_{\PQ}|$ are smooth.  A non--reduced component of a real member of
  $|\hol(2,2)|$ can only be of the form $2F$ with a real curve
  $F\in|\hol(0,1)|$. 
\end{lemma}

\pf
As usual, $\hol(k,l)$ denotes the locally free sheaf
$p_1^\ast\hol_{\bbfP^1}(k) \otimes p_2^\ast\hol_{\bbfP^1}(l)$ on
the smooth rational surface $\PQ$, where $p_i : \PQ \rightarrow \bbfP^1$ ($i
=1,2$) are the projections and $k, l$ are integers.
The Picard group $\Pic(\PQ)$ is free abelian of rank two with generators
$\hol(1,0)$ and $\hol(0,1)$. In the proof we shall use the well--known fact
that, if $k<0$ or $l<0$, then the linear system $|\hol(k,l)|$ is empty. 

Let $C\in|\hol(2,2)|$ be a real curve and $C_0\in |\hol(a,b)|$  an {\sl
irreducible\/} component (with {\sl reduced} scheme structure) of $C$. Let
$\lambda\geq 1$ be the multiplicity of $C_0$ in $C$, that is the largest
integer with $\lambda C_0\subset C$. Then we must have $0\leq \lambda a\leq 2$
and $0\leq \lambda b\leq 2$. The case $\lambda a = \lambda b = 2$ can only
occur if $\lambda C_0 = C$.

Assume first $\lambda C_0\ne C$, hence $\lambda a \leq 1$ or $\lambda b \leq
1$. If $\lambda\ge 2$ or $C_0$ singular, there exists a point $y\in C_0$ such
that any curve not contained in $C_0$ but containing $y$ has intersection
number at least two with $\lambda C_0$. If $F\in |\hol(0,1)|$ and
$G\in|\hol(1,0)|$ are the unique curves in these linear systems containing the
point $y$, we obtain (as $C_0$ is irreducible) $F.(\lambda C_0) = \lambda a\ge
2$ or $G.(\lambda C_0) = \lambda b\ge 2$. By the above inequalities, this means
$\lambda a =2$ and $0\leq\lambda b\leq 1$ or $0\leq\lambda a\leq 1$ and
$\lambda b =2$. If $\lambda = 1$, the curve $C_0$ is, by assumption,
irreducible, reduced, singular and a member of $|\hol(2,0)|$, $|\hol(2,1)|$,
$|\hol(0,2)|$, or $|\hol(1,2)|$. But these linear systems do not contain such a
curve. Hence, we must have $\lambda = 2$ and, therefore,  $C_0\in|\hol(0,1)|$
or $C_0\in|\hol(1,0)|$. In particular, $C_0$ is smooth.

If $C_0\in|\hol(1,0)|$, this curve is not real, since by definition of the real
structure on $\PQ$ this linear system does not contain real members. Hence, the
component $2C_0$ of $C$ is not real, which implies $2C_0 + 2\bar C_0 \subset
C$, since $C$ is real. But $\bar C_0 \in|\hol(1,0)|$ and so $2C_0 + 2\bar C_0
\in |\hol(4,0)|$. Such a curve can never be contained in $C\in|\hol(2,2)|$. So
we obtain $\lambda C_0 = 2F$ with some $F\in|\hol(0,1)|$.  
Again, since a curve of type $(0,4)$ can never be a component of $C$, the fibre
$F$ is necessarily real. This proves the lemma in the case $\lambda C_0 \ne C$.

Assume now $\lambda C_0=C$. Then $\lambda = 1$ or $\lambda = 2$. If
$\lambda=2$, we have, by reality of $C = 2C_0$, that  $C_0\in|\hol(1,1)|$ is a
real curve. Because $C_0$ is irreducible and $C_0.F =1$ for any
$F\in|\hol(0,1)|$, the curve $C_0$ would intersect each real fibre
$F\in|\hol(0,1)|$ at a real point. But on $\PQ$ real points do not exist. Hence
$\lambda=1$, which means $C = C_0$ is irreducible and real.

It remains to see that $C$ must be reducible if it is not smooth.  Let $x\in C$
be a singular point of $C$. Since $C$ is real and $\PQ$ does not contain real
points, $\bar{x}\ne x$ is also a singular point on $C$. If we embed $\PQ$ by
$|\hol(1,1)|$ as a smooth quadric into $\bbfP^3,$ we easily see that the linear
system of curves of type $(1,1)$ on $\PQ$ containing $x$ and $\bar{x}$ is
one--dimensional. It is cut out by the pencil of planes in $\bbfP^3$ containing
the line connecting $x$ and $\bar{x}$. Therefore, any point of $\PQ$ is
contained in such a curve.  The intersection number of $C$ with a curve of type
$(1,1)$ is four. Since $x$ and $\bar{x}$ are singular points on $C$, any curve
of type $(1,1)$ containing $x$ and $\bar{x}$ and a third point of $C$ must have
a common component with $C$. Therefore, $C$ cannot be irreducible and
reduced.\qed

\begin{defi}
  A reduced curve $C$ on a compact complex surface $S$ will be called a ``cycle
  of rational curves'', if the irreducible components $C_1,\dots, C_m$ of $C$
  are smooth rational curves with the following properties: (We use the
  convention $C_{m+1} = C_1$.)\\ 
  $m = 2$ and $C_1$ intersects $C_2$ transversally at two distinct points,
  $C_1.C_2 = 2$, or\\
  $m\geq 3$, $C_i.C_{i+1} = 1$ and $C_i\cap C_j \ne \emptyset$ implies
  $j\in\{i-1, i, i+1\}$.
\end{defi}

\begin{prop}\label{types}
Assume $\dim |\fund|\geq 1$ and let $S\in|\fund|$ be smooth and real. Then
there exists a blow-down $S\rightarrow \PQ$ and a connected real member
$C\in|-K_S|$, such that: 
\begin{itemize}
\item $\sigma$ is compatible with real structures, where we use the real
  structure of Lemma \ref{smocomp} on $\PQ$,
\item the composition $\mbox{pr}_2\circ \sigma$ of $\sigma$ with the second
  projection $\PQ\longrightarrow \bbfP^1$ is the morphism given by the linear
  system $|\widetilde{F}|$, where $\widetilde{F}\subset S$ is a real twistor
  fibre, 
\item the curve $C$ is reduced and
\item if $C$ is not smooth, it is a ``cycle of rational curves'' and its image
  $C'$ in $\PQ$ has one of the following structures:
\end{itemize}
\renewcommand{\labelenumi}{(\Roman{enumi})}
\begin{enumerate}
\item
\mbox{
\begin{minipage}[t]{10.5cm}
$C'$ has four components $C'=F+\bar{F}+G+\bar{G}$ where $F\in|\hol(0,1)|$ is a
non--real fibre and $G\in|\hol(1,0)|$ is a line.
\end{minipage}
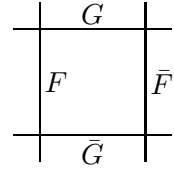
\begin{picture}(60,-60)(-15,-5)
\put(0,-40){\line(1,0){60}}  \put(30,-42){\makebox(0,0)[t]{$\bar{G}$}} 
\put(0,0){\line(1,0){60}}    \put(30,2){\makebox(0,0)[b]{$G$}}         
\put(10,10){\line(0,-1){60}} \put(12,-20){\makebox(0,0)[l]{$F$}}       
\put(50,10){\line(0,-1){60}} \put(52,-20){\makebox(0,0)[l]{$\bar{F}$}} 
\end{picture}\rule[-60pt]{0pt}{60pt}
}
\item
\begin{minipage}[t]{10.5cm}
$C'$ has two components $C'=F+C_0$ where $F\in|\hol(0,1)|$ is a real fibre and
$C_0\in|\hol(2,1)|$ is real, smooth and rational.
\end{minipage}
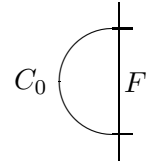
\begin{picture}(60,-60)(-15,-5)
\put(45,10){\line(0,-1){60}}  \put(47,-20){\makebox(0,0)[l]{$F$}}
\put(50,-20){\oval(55,40)[l]} \put(18,-20){\makebox(0,0)[r]{$C_0$}}
\end{picture}\rule[-60pt]{0pt}{60pt}
\item
\begin{minipage}[t]{13.8cm}
$C'$ has two distinct components $C'=A'+\bar{A}'$ where $A',\bar{A}'\in
|\hol(1,1)|.$\\
\mbox{\rm [By Corollary \ref{omit} this item can be omitted if
  $n=4$ and $\fdb$ is not nef!]}
\end{minipage}
\end{enumerate}
\renewcommand{\labelenumi}{(\roman{enumi})}
\end{prop}
\pf
As $\dim|-K_S| = \dim|\fund| - 1$, we have by assumption $|-K_S|\ne\emptyset$.
As seen above, we can choose a real blow--down map $\sigma:S\rightarrow\PQ$
such that $(\mbox{pr}_2\circ\sigma)^\ast\hol_{\bbfP^1}(1)
=\hol_S(\widetilde{F})$.  
If $|-K_S|$ contains a smooth member, we are done. Otherwise, take a reducible
real $C\in|-K_S|$ and let $C'\subset\PQ$ be the image of $C$. Since
$\sigma$ is a blow--up, $C'$ is a real member of $|\hol(2,2)|$. By Lemma
\ref{smocomp} the components of $C'$ are smooth and a multiple component can
only be the multiple $2F$ of a real fibre $F\in|\hol(0,1)|$.  But Lemma
\ref{noreal} shows 
that no point on such a real fibre is blown--up. Therefore, any other member of
$|2F|$ missing the $2n$ blown--up points, defines a divisor in $|-K_S|$.
Choosing, for example, a conjugate pair of appropriate fibres, we obtain a real
member in $|-K_S|$ whose image in $\PQ$ has only reduced components.

We assume for the rest of the proof that $C$ is chosen in this way.  If $C'$
would be irreducible, it would be smooth by Lemma \ref{smocomp}. In this case
$C$ is smooth, too.

Assume $C'$ is reducible.  A component of type $(1,a)$ with $a\in\{0,1,2\}$
cannot be real, since, otherwise, it would intersect real fibres at real
points. Therefore, such components appear in conjugate pairs, hence $a\le
1$. If $a=0$ we are in case (I).  If $a\ne 0$ then $a=1$ and $C'=A'+\bar{A}'$
with two distinct curves $A'$, $\bar{A}'$ of type $(1,1)$. This is case (III).

Assume now that there is no component of type $(1,a)$. Then we must have a
component $C_0$ of $C'$ which has type $(2,a)$. If $a=2$ it must be smooth and
we are done. Therefore, $a=1$, because $|\hol(2,0)|$ does not contain
irreducible reduced elements. Then we have $C'=C_0+F$ with $F\in|\hol(0,1)|$
and $C_0\in|\hol(2,1)|$. $F$ and $C_0$ must be real, since they have different
types. This is case (II) of our statement.

It remains to show that $C$ is a ``cycle of rational curves''. We have seen
this for the image $C'$ in $\PQ$.  Exceptional components of $C$ are always
rational. Furthermore, $C'$ has at most ordinary nodes as singularities. To
obtain $C$ from $C'$, at every step of blowing--up, we have to subtract the
exceptional locus from the total transform of $C'$. At every step we blow up
either a conjugate pair of singular points or of smooth points. We obtain a
curve which has again, at most, singularities of multiplicity two and is a
``cycle of rational curves''. So we obtain this property for $C$, too.\qed

Using this structure result and assuming that the fundamental linear system
$|\fund|$ is a pencil, we can show that the structure of its base locus is
closely related to the effective divisors of degree one on $Z$.
\begin{prop}\label{basel}
  Assume $\dim|\fund| = 1$ and denote by $C$ the base locus of the fundamental
  linear system.\\
  If $C$ is smooth, then $Z$ does not contain effective divisors of degree
  one.\\
  If $C$ is not smooth, then the number of effective divisors of degree one on
  $Z$  is equal to the number of components of $C$.
\end{prop}
\pf 
Let $S\in|\fund|$ be a smooth real member. Then $|-K_S| = \{C\}$ and by
Proposition \ref{types} $C$ is smooth or a ``cycle'' of smooth rational curves.
If $C$ is smooth, there does not exist an effective divisor $D$ of degree one,
because $D + \bar{D}\in|\fund|$ would produce a reducible member in
$|-K_S|$. Let now $C$ be singular, hence reducible.

The rest of the proof is an adaption of an idea of Pedersen and Poon
\cite[p.\ 700]{PP}.

Now let $\{P,\bar{P}\}$ be any pair of singular points on $C$. The image $C'$
of $C$ in $\PQ$ does not contain a real fibre, because, otherwise, by
Proposition \ref{types} and Lemma \ref{noreal} the linear system $|-K_S|$ would
be at least one--dimensional. Hence, the real twistor line $L_P$ containing $P$
and $\bar{P}$ is not contained in $S$. Hence, $L_P$ meets $S$ transversally at
$P$ and $\bar{P}$. If $Q$ is a point on $L_P$ distinct from $P$ and$\bar{P}$,
then there exists a divisor $S_0\in |\fund|$ containing $Q$.  Since $S_0$
contains also $C$ it contains three points of $L_P$. Hence, $L_P\subset S_0$.
Therefore, the real linear system of fundamental divisors containing $L_P$ is
non--empty. This implies that we can choose a real $S_0 \in |\fund|$ containing
$L_P$. 
Since $S_0$ contains also $C$ and $P$ is a singular point of $C$, the surface
$S_0$ contains three curves meeting at $P$, namely $L_P$ and two components
(call them $A$ and $B$) of $C$.  On the other hand, $L_P$ intersects $S$
precisely at $P$ and $\bar{P}$ as we have seen above. From $L_P.S=2$ we infer
that this is a transversal intersection. But $A$ and $B$ are contained in $S$
and are transversally there. We can conclude that the tangent space of $Z$ at
$P$ is generated by the tangent directions of $A, B$ and $L_P$ at $P$. Hence,
the real surface $S_0$ is singular at $P$. This implies that $S_0$ is singular
along $L_P$ (cf. \cite[p.141]{Hit2}) and by \cite[Lemma 2.1]{PP1} such a
divisor splits into the sum of two divisors of degree one.  Therefore, we have
at least as many pairs of conjugate divisors of degree one as we have pairs of
conjugate singular points on $C$. In other words, the number of distinct
divisors of degree one is at least equal to the number of components of $C$.

Let $D$ and $\bar{D}$ be a conjugate pair of divisors of degree one on $Z$.
Then $C\subset D\cup\bar{D}$. $D\cap \bar{D}$ is a real twistor line (cf.\
\cite[Prop. 2.1]{Ku}), and no component of $C$ is a real twistor line.
Hence, every component of $C$ lies on exactly one of the surfaces $D$ and
$\bar{D}$. By Lemma \ref{conn} $C\cap D$ is connected. The same is true for the
conjugate curve $C\cap \bar{D}$. Since $C$ is a cycle of rational curves,
$(C\cap D)\cap(C\cap\bar{D})$ consists of a conjugate pair $\{P,\bar{P}\}$ of
singularities of $C$. Since $D$ and $\bar{D}$ are of degree one, the real
twistor line $L_P$ containing $P$ and $\bar{P}$ must be contained in $D$ and in
$\bar{D}$. Therefore, $D\cap\bar{D}=L_P$.

Let $D'$ be an arbitrary divisor of degree one containing $L_P$. Then
$D'\cap\bar{D}'=L_P$ and without loss of generality we may assume $D'\cap
C=D\cap C$, since the decomposition of $C$ into two conjugate connected curves
is determined by $\{P, \bar{P}\}=L_P\cap C$. By Lemma \ref{pic} the restriction
map $\Pic{Z}\rightarrow\Pic{S}$ is injective. Since $(D'+\bar{D}')\cap S=C$ we
have $D'\cap S=D'\cap C$ and $D\cap S=D\cap C$. Hence, we have $\hol_Z(D)\cong
\hol_Z(D')$, which means that $D$ and $D'$ are linearly equivalent. If $D\ne
D'$, then $\dim|D|\ge 1$ and $Z$ would contain infinitely many divisors of
degree one and by Kurke \cite{Ku} and Poon \cite{Po2} it must be a
conic--bundle twistor space. But then we should have $\dim|\fund|=3$ in
contradiction to our assumption. Hence, $D=D'$ and we have exactly as many
divisors of degree one on $Z$ as $C$ has components.\qed

For technical reasons we state here the following lemma needed in Section
\ref{drei}. 
\begin{lemma}\label{normal}
  Let $S$ be a smooth complex surface and $C\subset S$ a reduced curve. Assume
  $C =\sum_{i=1}^{m} C_i$ is a ``cycle of rational curves'' as defined above.
  If $L\in\Pic(S)$ is a line bundle, we define $l_i := L.C_i$. Let $I_\pm :=
  \{i | \pm l_i>0\}$ and $C_\pm := \sum_{i\in I_\pm} C_i$. Let $\gamma$ denote
  the number of connected components of $C\setminus C_-$. Assume $|I_-|\geq 2$
  and each connected component of $C\setminus C_-$ contains a component of
  $C_+$. Then we have:
  \[h^0(C,L) = \sum_{i\in I_+} l_i - \gamma.\]
\end{lemma}
\pf 
Let $\eta:\tilde C = \sqcup_i C_i \rightarrow C$ be the normalization of $C$.
By $P_i$ we denote the intersection point of $C_i$ with $C_{i+1}$ $(1\leq i
\leq m)$. By assumption $m\geq 3$. Tensoring the exact sequence $0\rightarrow
\hol_C \rightarrow \eta_\ast \hol_{\tilde C} \rightarrow \oplus_i \bbfC_{P_i}
\rightarrow 0$ with $L$ yields the exact sequence $0 \rightarrow L\otimes
\hol_C \rightarrow \eta_\ast\eta^\ast(L\otimes\hol_C) \rightarrow \oplus_i
\bbfC_{P_i} \rightarrow 0$. Hence, $H^0(C, L\otimes\hol_C) \cong \ker(\oplus_i
H^0(C_i, L_i) \stackrel{\rho}{\rightarrow} \oplus_i \bbfC_{P_i})$. Here we
denote $L_i:= L\otimes\hol_{C_i} \cong \hol_{C_i}(l_i)$. Let $P'_i\in C_i$ and
$P''_i\in C_{i+1}$ be the two points on $\tilde C$ lying over $P_i$.

To describe $\rho$ we observe that the map $\eta$ gives isomorphisms $L_i(P'_i)
\stackrel{\sim}{\rightarrow} \bbfC_{P_i}$ and $L_{i+1}(P''_i)
\stackrel{\sim}{\rightarrow} \bbfC_{P_i}$. If $s_i\in H^0(C_i, L_i)$ is a
section, we denote by $s_i(P_i)$ the image of $s_i$ under the map $H^0(C_i,
L_i) \rightarrow L_i(P'_i) \stackrel{\sim}{\rightarrow} \bbfC_{P_i}$.
Similarly, $s_i(P_{i-1})$ is the image of $s_i$ under $H^0(C_i, L_i)
\rightarrow L_i(P''_{i-1}) \stackrel{\sim}{\rightarrow} \bbfC_{P_{i-1}}$. With
this notation we have:
\[ \rho(s_1,\dots, s_m) = (s_1(P_1) - s_2(P_1), s_2(P_2) - s_3(P_2),\dots,
s_m(P_m) - s_1(P_m)).\] Since $P'_i \ne P''_{i-1}$ on $C_i \cong \bbfP^1$, the
restriction of $\rho$ $H^0(C_i, L_i) \rightarrow
\bbfC_{P_i}\oplus\bbfC_{P_{i-1}}$ is surjective if and only if $l_i > 0$. If
$C_i + \cdots + C_{i+r}$ is a connected component of $C\setminus C_-$, then we
obtain by induction on $r$ that the restriction of $\rho$ $\oplus_{\mu = 0}^r
H^0(C_{i+\mu}, L_{i+\mu}) \rightarrow \oplus_{\mu = -1}^r \bbfC_{P_{i+\mu}}$ is
surjective. Because $H^0(C_i, L_i) = 0$ if and only if $l_i < 0$, we obtain
$\im(\rho) = \sum_{P_\mu\in C_0+C_+} \bbfC_{P_\mu}$. (Here, we denote $I_0 := I
\setminus (I_- \cup I_+)$ and $C_0 := \sum_{\nu\in I_0} C_\nu$.)
The number of points
$P_\mu\in C_0+C_+$ is equal to $|I_0| + |I_+| + \gamma$. Therefore, we obtain 
\[ \dim\ker(\rho) = \sum_i h^0(C_i, L_i) - (|I_0| + |I_+| +\gamma) = 
\sum_{l_i\geq 0} (l_i + 1) - |I_0| - |I_+| -\gamma = \sum_{l_i > 0} l_i -
\gamma.\] 
\qed
%
%
%
\section{The nef case}
\label{zwei}

For the rest of the paper we assume $n = 4$. Remember that $(\fund)^3 = 0$, 
$\chi(\fb{m}) = m+1$ and $h^0(\fdb) \geq 2$ in this case. Remember from
Mori's theory that a line bundle $L\in\Pic(Z)$ is called {\em nef}, if for each
irreducible curve $C\subset Z$ there holds $L.C\geq 0$.
\begin{thm}\label{nef}
The following properties are equivalent:
\begin{enumerate}
\item
$\fdb$ is nef;
\item
for all smooth and real $S\in|\fund|$ and all $C\in|-K_S|$, every component
$C_0$ of $C$ has the property $C_0.(-K_S)=0$;
\item
there exists a smooth and real $S\in|\fund|$ and a divisor $C\in|-K_S|$, such
that all components $C_0$ of $C$ have the property $C_0.(-K_S)=0$.
\end{enumerate}
If $\fdb$ is nef, then $a(Z)\le 2$ and $\dim|\fund|\leq 2$.\\ 
If $\fdb$ is nef and $\dim|\fund| = 2$, then $a(Z) = 2$, $|\fund|$ does not
have base points and for any smooth real $S\in|\fund|$ the pencil $|-K_S|$
contains a smooth real member.
\end{thm}

\pf
(i)$\Rightarrow$(ii):\\ Take any smooth real $S\in|\fund|$ and an arbitrary
curve $C\in|-K_S|$.  Since $C.(\fund)=0$ and $\fdb$ is nef we obtain (ii).

(ii)$\Rightarrow$(iii) is obvious.

(iii)$\Rightarrow$(i):\\ If $\fdb$ were not nef, then there would exist an
irreducible curve $C_0\subset Z$ with $C_0.(\fund)<0$. If $S\in|\fund|$ is
smooth and real, then $C_0\subset S$ and $C_0.(-K_S)=C_0.(\fund)<0$.
Therefore, $C_0$ is a component of any element of $|-K_S|$ in contradiction to
(iii).

Assume for the rest of the proof that $\fdb$ is nef. Let $S\in|\fund|$ be
smooth and real and $C\in|-K_S|$ a real member. If $C$ is smooth, we have shown
in \cite{CK} that $a(Z)\leq 2$ and $\dim|\fund|\leq2$. Assume $C$ is not
smooth. To compute the algebraic dimension consider the exact sequences
\[0 \rightarrow \fb{m-1} \rightarrow \fb{m} \rightarrow \can{m} \rightarrow 0\]
and 
\[0 \rightarrow \can{(m-1)} \rightarrow \can{m} \rightarrow N^{\otimes m}
\rightarrow 0\] 
with $N:=\canS\otimes\hol_C$. Since $\fdb$ is nef, $(-K_S).C_i = 0$ for any
component $C_i$ of $C$. But $C_i \cong \bbfP^1$ and so $N\otimes\hol_{C_i}
\cong \hol_{C_i}$. This does not imply in general $N \cong \hol_C$, because $C$
is a ``cycle'' of rational curves. But we obtain $h^0(N^{\otimes m}) = 1$ if
$N^{\otimes m} \cong \hol_C$ and $h^0(N^{\otimes m}) = 0$ if $N^{\otimes m}
\ncong \hol_C$. As in \cite{CK} this implies $a(Z)\leq 2$ and $h^0(\fdb) =
h^0(\hol_Z) + h^0(\canS) = 1 + h^0(\hol_S) + h^0(N) \leq 3$.

Assume now $\dim|\fund| = 2$. Hence, using the Riemann--Roch formula we have
$h^1(\fdb) = 1$. By Proposition \ref{tau} this implies $a(Z)\geq 2$. From the
above considerations we obtain $h^0(N) = 1$, hence $\canS\otimes\hol_C \cong N
\cong \hol_C$. (The same is true if $C$ is smooth, cf.\ \cite{CK}.) The exact
sequence $0 \rightarrow \hol_S \rightarrow \canS \rightarrow N \rightarrow 0$
and $h^1(\hol_S) = 0$ give a surjective restriction map $H^0(\canS)
\twoheadrightarrow H^0(\hol_C) \cong \bbfC$. Because $C\in|-K_S|$, this shows
that $|-K_S|$ does not have base points. Since $\dim|-K_S| = \dim|\fund| -1 =
1$, Bertini's Theorem \cite[I \S 1]{GH} states the existence of a smooth member
in $|-K_S|$. Hence, the generic divisor in $|-K_S|$ is smooth and so the
generic real member, too. On the other hand, we know from $h^1(\hol_Z) = 0$
that the restriction map $|\fund| \twoheadrightarrow |-K_S|$ is
surjective. From the freeness of  $|-K_S|$ we conclude that $|\fund|$ does
not have base points.\qed 

\begin{rem}\label{smnef}
  If there exists a smooth real $S\in|\fund|$ and a smooth curve $C\in|-K_S|$,
  then $\fdb$ is nef. This is clear from the theorem, because $C.(\fund) =
  (\fund)^3 = 0$.
\end{rem}

\begin{cor}\label{omit}
In Proposition \ref{types} we can omit case (III) if $\fdb$ is not nef.
\end{cor}

\pf
Assume $C'=A'+\bar{A}'$ as in the proof of Proposition \ref{types}.
$A'$ and $\bar{A}'$ intersect at a pair of conjugate points, say
$P$ and $\bar{P}$.

If $\sigma$ does not blow up $P$ and $\bar{P}$, then, by reality of the
blown--up set, on (the strict transforms of) $A'$ and $\bar{A}'$ exactly four
points are blown--up.  If we denote by $A$ and $\bar{A}$ the strict transforms
of $A'$ and $\bar{A}'$ in $S$, then we have $C=A+\bar{A}$ and
$A^2=\bar{A}^2=-2$.  Since $A$ and $\bar{A}$ are rational we obtain, by the
adjunction formula, $A.(-K_S)=\bar{A}.(-K_S)=0$. By Theorem \ref{nef} this
implies that $\fdb$ is nef.

If $\sigma$ blows up $P$ and $\bar{P}$, then we perform an elementary transform
to arrive at case (I) as follows. Let $\sigma_1:S^{(1)}\rightarrow\PQ$ be the
blow--up of 
$P$ and $\bar{P}$, then we have an induced real structure on $S^{(1)}$. Since
$A'$ intersects any fibre at exactly one point, $P$ and $\bar{P}$ lie on a
conjugate pair of fibres. The strict transforms in $S^{(1)}$ of these fibres
form a conjugate pair of disjoint $(-1)$--curves. Contracting them we obtain a
blow--down map $\sigma_1':S^{(1)}\rightarrow\PQ$ which is again compatible with
real structures.  If we denote by $E$ and $\bar{E}$ the exceptional curves of
$\sigma_1$, then the image of $C$ in $S^{(1)}$ is
$\sigma_1^\ast(A'+\bar{A}')-E-\bar{E} = A^{(1)}+\bar{A}^{(1)}+E+\bar{E}$. Here
$A^{(1)}$ and $\bar{A}^{(1)}$ are the strict transforms of $A'$ and $\bar{A}'$.
The morphism $\sigma_1'$ maps this curve onto a curve of type (I).\qed

%
%
\section{The non--nef case}
\label{drei}

Throughout this section we assume $\fdb$ to be not nef and $n=4$.\\
By Theorem \ref{nef}, Remark \ref{smnef} we know that in any smooth real
$S\in|\fund|$ the anticanonical system $|-K_S|$ contains only reducible
elements. By Proposition \ref{types} and Corollary \ref{omit} we can,
therefore, choose a real blow--down $S\rightarrow\PQ$ and a real reduced curve
$C\in|-K_S|$, whose image $C'$ in $\PQ$ is of type (I) or (II) as described
there. Observe that $C'$ has type (I) if and only if $C$ contains a real
irreducible component. 

\begin{prop}\label{notnefi}
  If there exists a real irreducible curve intersecting $|\fdb|$ negatively,
  then  $h^0(\fdb)=3$ and $a(Z) = 3$.\\
  There exists a unique irreducible curve $C_0$ with $C_0.(\fund) < 0$. This
  curve is real, smooth and rational and $C_0.(\fund)=-2$.  The base locus of
  $|\fund|$ is exactly $C_0$.  $Z$ does not contain divisors of degree one.
\end{prop}

\pf
Let $S\in|\fund|$ be smooth and real and choose $C\in|-K_S|$ and $\sigma:
S\rightarrow \PQ$ with the properties of Proposition \ref{types}. Because any
irreducible curve intersecting $\fdb$ negatively is contained in $C$, this
curve has a real component. Therefore, the image $C'$ of $C$ in $\PQ$ is of
type (II).
Let $C'=C_0'+F'$ be the decomposition of $C'$. By Lemma \ref{noreal} none of
the blown--up points lie on the real fibre $F'$. In particular, only smooth
points of $C_0'$ are blown--up. Hence, $C=C_0+F$ where $C_0$ and $F$ are the
strict transforms of $C_0'$ and $F'$ respectively.  Therefore, the eight
blown--up points lie on $C_0'$ which implies $C_0^2 = {C'}_0^2-8=-4$.

By adjunction formula we obtain $-2=C_0.(-K_S)=C_0.(\fund)$. Hence,
$|-K_S|=C_0+|F|$ and we obtain: $\dim|-K_S|=1$ and $C_0$ is the base locus of
$|-K_S|$. Since $h^1(\hol_Z)=0$ the restriction map
$H^0(\fdb)\twoheadrightarrow H^0(\canS)$ is surjective. Hence, the linear
system $|\fund|$ has dimension two and its base locus is precisely $C_0$ (with
multiplicity one). $C_0$ is the unique irreducible curve in $Z$ having negative
intersection number with $\fund$, since any other such curve should be
contained in the base locus of $|\fund|$.

If $Z$ contains a divisor $D$ of degree one, then $D+\bar{D}\in|\fund|$. If
$D_0$ ($\bar{D}_0$ respectively) denotes the restriction of $D$ to $S$, then
$D_0+\bar{D}_0\in|-K_S|=C_0+F$.  In the proof of Lemma \ref{noreal} we have
seen that the real elements of $|F|$ are irreducible. Therefore, any real
element of $|-K_S|$ consists of two distinct real irreducible curves with
multiplicity one. This shows that $|-K_S|$ cannot contain a member of the form
$D_0+\bar{D}_0$. 

It remains to show that the {\bf algebraic dimension} of Z must be three in
this case.

Let $\sigma:\tilde{Z}\rightarrow Z$ be the blow--up of the smooth rational
curve $C_0$. By $E\subset \tilde{Z}$ we denote the exceptional divisor. Then we
obtain a morphism $\pi:\tilde{Z}\rightarrow\bbfP^2$ defined by the linear
system $|\fund|$ such that $\pi^\ast\hol(1)\cong
\sigma^\ast\fdb\otimes\hol_{\tilde{Z}}(-E)$.
Since the restriction map $|\fund|\rightarrow|-K_S|$ is surjective, the
restriction $\pi_{|S}$ is given by the linear system $|-K_S|=C_0+|F|$. This
means that $\pi$ exibits $S$ as the blow--up of a ruled surface and $\pi(S)$ is
a line in $\bbfP^2$. Since $\pi(\tilde{Z})$ is not contained in a linear
subspace, $\pi$ must be surjective. If we equip $\bbfP^2$ with the usual real
structure, $\pi$ becomes compatible with real structures since the linear
system $|\fund|$ and the blown--up curve $C_0$ are real.

Since $Z$ does not contain divisors of degree one, any real fundamental divisor
$S$ is irreducible and, therefore, smooth. By $\tilde{S}\subset\tilde{Z}$ we
denote the strict transform of $S\in |\fund|$. Since $C_0$ is a smooth curve in
a smooth surface, $\sigma:\tilde{S}\rightarrow S$ is an isomorphism.
Furthermore, $E\cap\tilde{S}$ will be mapped isomorphically onto $C_0\subset
S$.  Since $F.C_0=2$ and the restriction of $\pi$ onto $\tilde{S}$ is the map
defined by the linear system $|F|$, the restriction of $\pi$ exibits
$E\cap\tilde{S}$ as a double covering over $\pi(S)\cong\bbfP^1$.  Since real
lines cover $\bbfP^2$ the morphism $\pi:E\rightarrow\bbfP^2$ does not contract
curves and is of degree two.

Since generic fibres of $\pi$ are smooth rational curves, the line bundle
$\hol_{\tilde{Z}}(E)$ restricts to $\hol_{\bbfP^1}(2)$ on such fibres.  Hence,
after replacing (if necessary) $\bbfP^2$ by the open dense set $U$ of points
having smooth fibre, the adjunction morphism $\pi^\ast\pi_\ast
\hol_{\tilde{Z}}(E) \rightarrow \hol_{\tilde{Z}}(E)$ is surjective. This
defines a $U$--morphism $\Phi:\tilde{Z} \rightarrow
\bbfP(\pi_\ast\hol_{\tilde{Z}}(E))$, where  $\pi_\ast\hol_{\tilde{Z}}(E)$ is a
locally free sheaf of rank three.  The restriction of $\Phi$ to smooth fibres
coincides with the Veronese embedding $\bbfP^1\hookrightarrow\bbfP^2$ of degree
two. Therefore, the image of $\Phi$ is a three--dimensional subvariety of the
$\bbfP^2$--bundle $\bbfP(\pi_\ast\hol_{\tilde{Z}}(E))\rightarrow U$.  Hence,
$\tilde{Z}$ is bimeromorphically equivalent to a quasiprojective variety and
has, therefore, algebraic dimension three.\qed

For the rest of this section we assume that there does not exist a {\em real}
irreducible curve contained in the base locus of $|\fund|$. We keep the
assumptions $n=4$ and $\fdb$ is not nef. In this situation, we obtain:
\begin{lemma}\label{sub}
  \begin{enumerate}
  \item[(a)] If $A\subset Z$  is an irreducible curve, then $A.(\fund) \geq
    -2$. 
  \item[(b)] If $A\subset Z$ is an irreducible curve with $A.(\fund) <0$, then
    there 
    exists at least a one--parameter family of real smooth divisors
    $S\in|\fund|$, containing a curve $C\in|-K_S|$ and possessing a
    birational morphism $\sigma: S\longrightarrow \PQ$ as in Proposition
    \ref{types}, such that moreover:\\
    $A$ and $\overline{A}$ are components of $C$ and for twistor fibres
    $F\subset S$ we have $F.A = F.\overline{A} = 1$.\\
    In particular, the image $A'$ of $A$ in $\PQ$ is a ``line'', that means
    $A'\in|\hol(1,0)|$. 
  \end{enumerate}
\end{lemma}

\pf 
Our assumptions imply that $C'$ is a curve of type (I) in Proposition
\ref{types}. The components of $C'$ are curves in $\PQ$ with self--intersection
number zero. They are not real. Hence, after a succession of four blow--ups of
a 
conjugate pair of points, each component $A$ of $C$ fulfills $A^2\geq-4$ in
$S$. The adjunction formula, together with the rationality of $A$, implies
$A.(\fund) = A.(-K_S) = A^2 + 2 \geq -2$. Because a curve $A$ with $A.(\fund) <
0$ must be a component of $C$, the assertion (a) is shown. 

Let now $A\subset Z$ be an irreducible curve with $A.(\fund)<0$. Then we have
$A\subset C$. Let $x\in A\subset C$ be a smooth point of $C$ and $x\in F
\subset Z$ a twistor fibre. Since $|\fund|$ is at least a pencil, there exists
a divisor $S\in|\fund|$ containing a given point $y\in
F\setminus\{x,\overline{x}\}$. Because $F.S = 2$ and $S\cap F \supset \{y, x,
\overline{x}\}$ the twistor fibre $F$ is contained in $S$. So the real
subsystem $|\fund|_F\subset |\fund|$ of divisors containing $F$ is not
empty. Hence, we can choose a real smooth $S\in|\fund|$ containing $F$. By
construction, we have $F.A = F.\overline{A} \geq 1$. But $F.B\geq 0$ for any
curve $B\subset S$ together with $F.(-K_S)=2$ implies $F.A = F.\overline{A} = 
1$. Because $S$ contains only a real one--parameter family of real twistor
lines, the intersection points with real twistor fibres form only a real
one--dimensional subset of points $z$ on $A$. Therefore, we obtain at least a
one--parameter family of such surfaces $S$. Proposition \ref{types} implies now
the claim. 
\qed

\begin{prop}\label{notnefii}
Assume the existence of an irreducible (non--real) curve $A\subset Z$ with
  $A.(\fund) = -2$. Then:   $h^0(\fdb)=4$ and $a(Z) = 3$.\\
  The curves $A$ and $\bar{A}$ are disjoint smooth and rational.  $A$ and 
  $\bar{A}$ are the unique irreducible reduced curves having negative
  intersection number with $\fund$.  The base locus of $|\fund|$ is exactly the
  union of $A$ and $\bar{A}$.  $Z$ contains infinitely many divisors of degree
  one and is one of the twistor spaces studied by LeBrun
  \cite{LeB} and Kurke \cite{Ku}.
\end{prop}

\pf 
We choose $S\in|\fund|$ as in Lemma \ref{sub}(b). Then we have $C\in|-K_S|$
containing $A$ and $\overline{A}$ as smooth rational components which are
mapped to ``lines'' $A'$ and $\overline{A}'\in|\hol(1,0)|$ in $\PQ$. We have by
the adjunction formula $A^2 = A.(-K_S) -2 = -4$, which implies that the eight
blown--up points lie on $A'+\overline{A}'$ (or the strict transforms after
partial blow--ups). Hence, any member of $A'+\overline{A}'+|\hol(0,2)|$
contains the eight blown--up points and defines, therefore, a divisor in
$|-K_S|$. On the other hand, any curve in $|-K_S|$ is mapped onto a curve of
type $(2,2)$ on $\PQ$ containing the blown--up points. Such a curve must
contain $A'$ and $\overline{A}'$ since the intersection number with $A$ is two,
but four of the blown--up points lie on $A$. Hence, the image of $|-K_S|$ is
precisely the two--dimensional system $A'+\overline{A}'+|\hol(0,2)|$. So we
have $\dim|-K_S|=2$. Using the exact sequence $0 \rightarrow \hol_Z \rightarrow
\fdb \rightarrow \canS \rightarrow 0$, this implies $h^0(\fdb)=4$. Furthermore,
we see that $A+\bar{A}$ is in the base locus of $|-K_S|$ which coincides with
the base locus of $|\fund|$.

To see that $Z$ contains infinitely many divisors of degree one we modify an
idea of Pedersen, Poon \cite[p.\ 700]{PP}:

The above description of the image of $|-K_S|$ in $\PQ$ shows that there exist
infinitely many real curves $C\in|-K_S|$ whose image $C'$ in $\PQ$ is of type
(I) and the singular points of $C'$ are not blown--up. Then $C$ is the strict
transform of  such a curve  and $C'$ consists of four irreducible components. 

Let $P$ and $\bar{P}$ denote a conjugate pair of singular points of $C$. Since
$A$ and $\bar{A}$ are disjoint, both points $P$ and $\bar{P}$ are contained in
$A+\bar{A}$, hence in the base locus of $|\fund|$. On the twistor space $Z$
there exists exactly one real twistor line $L_P$ connecting $P$ with $\bar{P}$.
Let $Q\in L_P$ be a point different from $P$ and $\bar{P}$. The linear system
$|\fund|_Q$ of all fundamental divisors containing $Q$ has at least dimension
$\dim|\fund|-1=\dim|-K_S|=2$. Since fundamental divisors have degree two and
any member of $|\fund|_Q$ contains $P,\,\bar{P}$ and $Q$ it also contains
$L_P$. Hence, $|\fund|_Q$ coincides with the {\sl real} linear subsystem of
divisors in $|\fund|$ containing $L_P$.

Choose now a point $R\in C$ which is not on $A+\bar{A}$. Then the real linear
system $|\fund|_{L_P,R,\bar{R}}$ is non--empty.
Let $C$ be decomposed as $A+\bar{A}+B+\bar{B}$. Then, by our choice of $C$,
$B.(\fund)=2$ and $A.(\fund)=-2$. These four curves intersect as indicated in
the following picture:
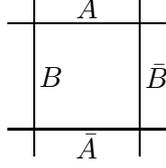
\begin{picture}(90,-90)(-30,0)
\put(0,-40){\line(1,0){60}}  \put(30,-42){\makebox(0,0)[t]{$\bar{A}$}} 
\put(0,0){\line(1,0){60}}    \put(30,2){\makebox(0,0)[b]{$A$}}         
\put(10,10){\line(0,-1){60}} \put(12,-20){\makebox(0,0)[l]{$B$}}       
\put(50,10){\line(0,-1){60}} \put(52,-20){\makebox(0,0)[l]{$\bar{B}$}} 
\end{picture}
\rule[-50pt]{0pt}{60pt}
\vspace*{12pt}

Any real member $S_0$ of $|\fund|_{L_P,R,\bar{R}}$ contains three distinct
points of $B$, namely $R$ and the intersection of $B$ with $A+\bar{A}$. Hence,
$B\subset S_0$. But this means that $S_0$ contains three curves, say $A, B$ and
$L_P$, which meet at $P\in S_0$.  As in the proof of Proposition \ref{basel} we
obtain that $S_0$ is reducible, hence splits into the sum of two divisors of
degree one.  Since the intersection of a conjugate pair of divisors of degree
one is a twistor line, we obtain in this way infinitely many divisors of degree
one on $Z$.

By the result of Kurke--Poon \cite{Ku}, \cite{Po2} we obtain that $Z$ is a
LeBrun twistor space and $A+\bar{A}$ is precisely the base locus of
$|\fund|$. Hence, no other curve can have negative intersection number with
$(\fund)$. This proves that really all properties of the Propositon are
fulfilled.\qed

\begin{prop}\label{notnefiii}
  Assume $A.(\fund)\geq-1$ for all irreducible curves $A\subset Z$. Then:
  $h^0(\fdb)=2$ and $a(Z) = 3$.\\
  There exists a conjugate pair of irreducible curves $A$ and $\bar{A}$ which
  are smooth and rational and $A.(\fund)=\bar{A}.(\fund)=-1$.  The base locus
  $C$ of $|\fund|$ consists of a cycle of an even number of rational curves.
  The number of distinct divisors of degree one on $Z$ is equal to the number
  of components of $C$.
\end{prop}

\pf
By our assumption we obtain the existence of an irreducible curve $A\subset Z$
with $A.(\fund) = -1$. Now we choose $S\in|\fund|$ real and smooth and
$C\in|-K_S|$ as in Lemma \ref{sub}(b). The curve $C\in|-K_S|$ has $A$ and
$\overline{A}$ as components and the images $A'$ and $\overline{A}'$ of $A$ and
$\overline{A}$ in $\PQ$ are members of $|\hol(1,0)|$. But $A^2 = A.(-K_S) - 2 =
-3$ implies that exactly one pair of blown--up points does not lie on $A' +
\overline{A}'$. This implies that the two components of $C'$, which are members
of $|\hol(0,1)|$, are not movable. Hence, $|-K_S| = \{C\}$ and, as above,
$h^0(\fdb)=2$. Because $C'$ is of type (I), the curve $C$ consists of 2, 3, 4,
5 or 6 pairs of conjugate rational curves.

By Proposition \ref{basel} we have: the number of components of $C$ is equal
to the number of effective divisors of degree one.

It remains to be seen that the {\bf algebraic dimension} is three.
Since $|\fund|$ is a pencil we study in more detail the linear system $|-K|$ on
$Z$.

We need to investigate the structure of $C$ before we can collect more
information on the linear system $|-2K_S|$.

We know that the blow--up $\sigma:S\rightarrow S^{(0)} := \PQ$ factors through
a 
succession of four blow--ups $S=S^{(4)} \rightarrow S^{(3)} \rightarrow S^{(2)}
\rightarrow S^{(1)} \rightarrow S^{(0)}$ such that at each step a conjugate
pair of points is blown--up. The image of $C$ in $S^{(i)}$ will be denoted by
$C^{(i)}$.  The blown--up points in $S^{(i)}$ should lie on $C^{(i)}$. If they
are smooth points of $C^{(i)}$ then $C^{(i+1)} \stackrel{\sim}{\rightarrow}
C^{(i)}$.  If we blow up a conjugate pair of singular points of $C^{(i)}$, the
curve $C^{(i+1)}$ has two components more than $C^{(i)}$. By assumption,
$C^{(0)} = C' \subset \PQ$ is of type (I). Each $C^{(i)}$ is a ``cycle of
rational curves''. We can choose the factorization of $\sigma$ in such a way
that at the first $k$ steps, only singular points of $C^{(i)}$ are blown--up
and at the last $4-k$ steps, only smooth points of $C^{(i)}$ are blown--up.
Then $C$ will have $2(2+k)$ components, where $0\leq k \leq 4$.  If we would
have a component $A$ of $C$ with $A^2 =0$, then the image 
$A^{(0)}$ of $A$ in $S^{(0)}$ would be a component of $C^{(0)}$ and none of the
blown--up points would lie on $A^{(0)}$. But then four of the blown--up points
must lie on a line or on a fibre in $S^{(0)}$, which implies that $C$ has a
component $B$ with $B^2 = -4$. This was excluded by assumption.  Therefore, for
any component $A$ of $C$ we have  $-1 \geq A^2 \geq -3$.
Since $A$ is a smooth rational curve, this means $A.(-K_S) \in \{-1, 0, +1\}$.

By assumption $C$ is reduced. Let $C = \sum\nolimits_{\nu = 1}^m C_\nu$ be the
decomposition of $C$ into irreducible components. By Proposition \ref{types} we
have $C_\nu \cong \bbfP^1$, $C_\nu.C_{\nu+1} = 1$ and $C_\nu$ intersects only
$C_{\nu-1}$ and $C_{\nu+1}$. (This means ``$C$ is a cycle of rational curves''.
For convenience, we use cyclic subscripts, that is $C_\nu = C_{\nu+m}$.) For
$\varepsilon\in\{-1, 0, +1\}$ we define $I_\varepsilon := \{\nu | C_\nu.(-K_S)
= \varepsilon\}$ and $C_\varepsilon := \sum\nolimits_{\nu\in I_\varepsilon}
C_\nu$. In this way we split $C$ into three parts $C = C_- + C_0 + C_+$. As
$C.(-K_S) = 0$ we have $|I_-| = |I_+|$.  The assumption that $\fdb$ is not nef
implies $I_- \ne \emptyset$. The curve $C$ has no real component, hence $|I_-|
= |I_+|\geq 2$.

{\bf Claim:} Any two components of $C_+$ are disjoint.\\ 
Let $C_\alpha$ and $C_\beta$ be two
distinct components of $C_+$. Then $C_\alpha^2 = C_\beta^2 = -1$. If $C_\alpha$
and $C_\beta$ are both 
contracted to a point on $S^{(0)}$, they are obviously disjoint.
If $C_\alpha$ and $C_\beta$ are mapped to curves in $S^{(0)}$, by our choice of
$S$ both must be members of $|\hol(0,1)|$, because $A$ and $\overline{A}$ are
components of $C_-$. Finally, we have to exclude the case  where $C_\alpha$ is
mapped to a curve $C'_\alpha\in|\hol(0,1)|$ in $\PQ$ and $C_\beta$ is
contracted to a point $P\in C'_\alpha$. This implies $C_\alpha \ne
\overline{C}_\beta$. Since $C_\beta$ is a component of the anticanonical
divisor 
$C\subset S$, the point $P$ must be singular on $C'$. Thus, we can take
$S^{(1)} \longrightarrow S^{(0)} = \PQ$ to be the blow--up of $P$ and
$\overline{P}$. The curve $C^{(1)}$ consists then of six $(-1)$--curves among
which we find the images of $C_\alpha, C_\beta, \overline{C}_\alpha$ and
$\overline{C}_\beta$. Because those are $(-1)$--curves on $S$, the remaining
six blown--up points must lie on one pair of $(-1)$--curves giving 
rise to $(-4)$--curves in $C$ contradicting our assumption.

Thus, the claim is proved and $C_+$ is the disjoint union of an even number of
smooth rational $(-1)$--curves.

Since $C$ is a cycle of rational curves, the curve $C\setminus C_+ = C_- + C_0$
has the same number of connected components as $C_+$.

We claim that each connected component of $C_- + C_0$ contains exactly one
component of $C_-$. Since $C_-$ and $C_+$ have the same number of components,
this is equivalent to the statement that each connected component of $C_- +
C_0$ contains at most one component of $C_-$.\\ Assume the contary, that is
there is a connected component of $C_- + C_0$ containing two irreducible
components of $C_-$. Then these two components of $C_-$ are not conjugate to
each other. (Two conjugate components of $C$ are ``opposite'' in the cycle
$C$.) Therefore, $C_-$ has at least four components and so $C_+$. Thus $C_- +
C_0$ has at least four connected components, hence at least six irreducible
ones. Therefore, $C$ contains at least ten irreducible components, that is
$k\geq 3$. Therefore, the image $C^{(3)}$ of $C$ in $S^{(3)}$ consists of a
cycle of ten rational curves with self--intersection numbers $-2, -1, -3, -1,
-2, -2, -1, -3, -1, -2$ (in this order). By assumption, at the last step of
blow--up, no point on a $(-3)$--curve is blown--up. To obtain a second pair of
$(-3)$--curves, we have to blow up points on a conjugate pair of
$(-2)$--curves. 
Such curves are not neighbours in the cycle $C^{(3)}$, they are opposite to
each other. Thus, one easily sees that after the last step of blow--up, between
two $(-3)$--curves on the cycle $C$ we always have a $(-1)$--curve. This means
that no connected component of $C\setminus C_+$ contains two irreducible
components of $C_-$, as claimed.

We can now compute the dimension of $|-K|$ and $|-2K_S|$.\\ The Riemann--Roch
formula and $h^0(\fdb) = 2$ imply $h^1(\fdb) = 0$. Hence, the exact sequence $0
\rightarrow \fdb \rightarrow K^{-1} \rightarrow \can{2} \rightarrow 0$ gives
$h^0(K^{-1}) = h^0(\fdb) + h^0(\can{2}) = 2 + h^0(\can{2})$ and a surjective
restriction map $|-K| \twoheadrightarrow |-2K_S|$. With $N:=\canS \otimes
\hol_C$ we obtain an exact sequence $0 \rightarrow \can{1} \rightarrow \can{2}
\rightarrow N^{\otimes 2} \rightarrow 0$. Since $h^1(\canS) = h^1(\fdb) = 0$,
this sequence yields $h^0(\can{2}) = h^0(\canS) + h^0(N^{\otimes 2}) = 1 +
h^0(N^{\otimes 2})$.  We can apply Lemma \ref{normal} with $L=\can{2}$, because
we have seen that the components of $C_+$ are disjoint to each other and that
each connected component of $C\setminus C_+$ contains exactly one irreducible
component of $C_-$. We obtain $h^0(N^{\otimes 2}) = h^0(C,L) = \sum_{\nu\in
  I_+} C_\nu.(-2K_S) - |I_-| = 2|I_+| - |I_-| = |I_+|$. Hence, $\dim |-2K_S| =
|I_+| \geq 2$ and $\dim |-K| = 2 + \dim |-2K_S| \geq 4$.

Next we study the base locus of $|-2K_S|$.\\ Since any component $A$ of $C_-$
fulfills $A.(-K_S) = -1$, $C_-$ is in the base locus of $|-2K_S|$. Now let $A$
be an irreducible component of $C_0$. Since, as we have shown above, any
connected component of $C_0 + C_-$ contains a curve $B\subset C_-$ there exists
a finite chain of components $A_1,\dots,A_r = A$ of $C_0$ with $B\cap A_1 \ne
\emptyset$ and $A_i \cap A_{i+1} \ne \emptyset\quad (1\leq i < r)$. But for a
component $A_i$ of $C_0$ we have $A_i.(-2K_S) = 0$ which implies: if $A_i$
intersects the base locus of $|-2K_S|$, it must be contained in this base
locus. Hence, by induction on $i$, we obtain that $A_i$ is contained in the
base locus of $|-2K_S|$ for all $1\leq i \leq r$. This shows that $C_0 + C_-$
is contained in the base locus of the linear system $|-2K_S|$.  Therefore, we
have $|-2K_S| = C_0 + C_- + |-K_S + C_+|$ and the map defined by $|-2K_S|$
coincides with the map given by $|-K_S + C_+|$.

We have seen above that $C_+$ is a disjoint union of smooth rational
$(-1)$--curves. We can, therefore, contract $C_+$ to obtain a smooth surface
$S'$. If $\sigma':S\rightarrow S'$ is this contraction, we have
$\canS\otimes\hol_S(C_+) \cong \sigma^{'\ast}(K^{-1}_{S'})$. Hence, $-K_{S'}$
is nef if and only if $-K_S + C_+$ is nef.

First we deal with the case where $C_+ - K_S$ is nef. In this case, the generic
member of the moving part of $|-2K_S|$ is irreducible. This can be seen as
follows: 

Observe first that $S'$ can be blown down to $\bbfP^2$.
This follows from \cite[Prop. 3, p. 48]{Dem} because $-K_{S'}$ is nef and
$K^2_{S'} = |I_+| > 0$. 

Because $C_+$ has 2, 4 or 6 components, $S'\rightarrow\bbfP^2$ is a blow--up of
7, 5 or 3 points. Therefore, we can apply a theorem of Demazure \cite[p.\ 39
and p.\ 55]{Dem} stating that $|-K_{S'}|$ contains a smooth irreducible member
and is base point free if $|-K_{S'}|$ is nef. Hence, there exists a smooth
irreducible curve in $|-K_{S'}|$ avoiding the blown--up points. Its preimage in
$S$ is a smooth irreducible member of $|\sigma^\ast K^{-1}_{S'}| = |C_+ - K_S|
= |-2K_S - C_0 -C_-|$ which is, therefore, the moving part of $|-2K_S|$.

Since $\dim |-2K_S| \geq 2$ and the generic member of the moving part of
$|-2K_S|$ is a smooth irreducible curve, the image of the map defined by
$|-2K_S|$ has dimension two. We have seen before that the restriction
$|-K|\twoheadrightarrow |-2K_S|$ is surjective which implies that the map
$\Phi_{|-K|}$ given by $|-K|$ on $Z$ coincides, after restriction to $S$, with
the map given by $|-2K_S|$. Since $2S\in|-K|$, the $\Phi_{|-K|}$--image of $S$
is contained in a hyperplane. But the $\Phi_{|-K|}$--image of $Z$ cannot be
contained in a hyperplane. Therefore, the image of $\Phi_{|-K|}$ has dimension
three. This implies $a(Z)=3$.

We are left with the case where $C_+ - K_S$ is not nef.\\
In this case we shall see that $\Phi_{|-K|}$ has only two--dimensional image
but equips $Z$ with a conic--bundle structure. Under the assumption that $C_+ -
K_S$ is not nef we study the structure of $C$.

Let $A$ be an irreducible curve in $S$ with $A.(C_+ - K_S) < 0$. If
$A\nsubseteq C$, then $A.C_+\geq 0$. But the base locus of $|-K_S|$ is
contained in $C$ and this implies $A.(-K_S)\geq 0$. Hence, we have necessarily
$A\subseteq C$.  If $A\subseteq C_+$, then $A.(-K_S) = 1$ and $A.C_+ \geq
A^2$, hence, $A.(C_+ - K_S) \geq A^2 + 1 = 0$.  If
$A\subseteq C_0$, then $A.(-K_S) = 0$ and $A.C_+\geq 0$, hence, $A.(C_+ - K_S)
\geq 0$.  If, finally, $A\subseteq C_-$, then $A.(C_+ - K_S) = A.C_+ - 1 \geq
-1$.  So, we obtain: the irreducible curves $A\subset S$ with $A.(C_+ - K_S) <
0$ are exactly those components of $C_-$ which are disjoint to $C_+$. They
fulfill $A.(C_+ - K_S) = -1$. 

As we have seen above, each connected component of $C\setminus C_+$ contains
exactly one irreducible component of $C_-$. Hence, if $A$ is a component of
$C_-$ which does not meet $C_+$, then its connected component should contain
at least two curves from $C_0$. Thus, using reality, we see that $C$ has at
least eight components. This means, using the convention introduced above, the
image $C^{(2)}$ of $C$ in $S^{(2)}$ (the surface
obtained after two steps of blow--up) consists of eight curves whose
self--intersection numbers are alternately $-1$ and $-2$:
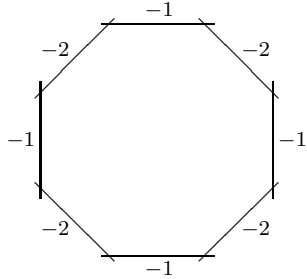
\vspace{2mm}
\centerline{
\begin{picture}(90,90)(0,-90)
\put(2,-32){\line(1,1){30}}  \put(15,-17){\makebox(0,0)[br]{$\scriptstyle -2$}}
\put(4,-70){\line(0,1){44}}  
\put(2,-48){\makebox(0,0)[r]{$\scriptstyle -1$}}     
\put(2,-64){\line(1,-1){30}} \put(15,-79){\makebox(0,0)[tr]{$\scriptstyle -2$}}
\put(27,-92){\line(1,0){43}} 
\put(49,-94){\makebox(0,0)[t]{$\scriptstyle -1$}}    
\put(64,-94){\line(1,1){30}} \put(80,-79){\makebox(0,0)[tl]{$\scriptstyle -2$}}
\put(92,-70){\line(0,1){44}} 
\put(94,-48){\makebox(0,0)[l]{$\scriptstyle -1$}}    
\put(94,-32){\line(-1,1){30}}\put(80,-17){\makebox(0,0)[bl]{$\scriptstyle -2$}}
\put(27,-4){\line(1,0){43}}  
\put(49,-2){\makebox(0,0)[b]{$\scriptstyle -1$}}     
\end{picture}
}
\vspace{1pt}

One now easily sees: if we were to blow up a pair of singular points on
$C^{(2)}$, in the resulting curve $C^{(3)}$ any $(-2)$--curve would meet a
$(-1)$--curve and any $(-3)$--curve would meet two $(-1)$--curves. Therfore,
after the last step of blow--up, no $(-3)$--curve is disjoint to all
$(-1)$--curves on $C$. Thus, we can only blow up smooth points in the last two
steps. If we were to blow up a conjugate pair of points on $(-2)$--curves, the
resulting $(-3)$--curve would intersect two $(-1)$--curves. Then, again, in $C$
there would be no $(-3)$--curve disjoint to all $(-1)$--curves. So we conclude
that the last two pairs 
of conjugate blown--up points cannot lie on $(-2)$--components of $C^{(2)}$. If
each of the four $(-1)$--curves in $C^{(2)}$ contains one of the blown--up
points, then any component of $C$ has zero intersection number with $-K_S$.
But, then by Theorem \ref{nef}, $\fdb$ would be nef which contradicts our
general assumption. Thus, the four blown--up points lie on a pair of conjugate
$(-1)$--curves. The structure of $C$ is, therefore, the following:

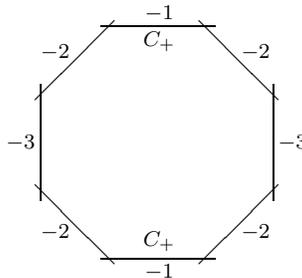
\centerline{
\begin{picture}(90,90)(0,-90)
\put(2,-32){\line(1,1){30}}  \put(15,-17){\makebox(0,0)[br]{$\scriptstyle -2$}}
\put(4,-70){\line(0,1){44}}  
\put(2,-48){\makebox(0,0)[r]{$\scriptstyle -3$}}      
\put(2,-64){\line(1,-1){30}} \put(15,-79){\makebox(0,0)[tr]{$\scriptstyle -2$}}
\put(27,-92){\line(1,0){43}} 
\put(49,-94){\makebox(0,0)[t]{$\scriptstyle -1$}}     
\put(49,-89){\makebox(0,0)[b]{$\scriptstyle C_+$}}    
\put(64,-94){\line(1,1){30}} \put(80,-79){\makebox(0,0)[tl]{$\scriptstyle -2$}}
\put(92,-70){\line(0,1){44}} 
\put(94,-48){\makebox(0,0)[l]{$\scriptstyle -3$}}     
\put(94,-32){\line(-1,1){30}}\put(80,-17){\makebox(0,0)[bl]{$\scriptstyle -2$}}
\put(27,-4){\line(1,0){43}}  
\put(49,-2){\makebox(0,0)[b]{$\scriptstyle -1$}}      
\put(49,-6){\makebox(0,0)[t]{$\scriptstyle C_+$}}     
\end{picture}
}
\vspace{11pt}

In particular, we obtain: $\dim |-2K_S| = |I_+| = 2$ and $\dim |-K| = 2 +
\dim|-2K_S| = 4$.
Furthermore, since both components of $C_-$ have negative intersection number
with $C_+ - K_S$, the curve $C_-$ is contained in the base locus of $|C_+ -
K_S|$. This means $|-2K_S| = |2C_+ + C_0| + C_0 + 2C_-$.
By our choice of $S$, the two components $A$ and $\overline{A}$ of $C_-$ are
mapped onto lines $A'$ and $\overline{A}'$ on $\PQ$. The above analysis of $C$
shows that we can decompose $\sigma : S \longrightarrow \PQ$ into the
following steps: 

First we blow up a conjugate pair of singular points on the
curve $C'$ (which is of type (I)). This produces precisely two singular fibres
of the ruling (whose general fibre is the image of a twistor fibre). In the
second step we blow up the two singular points of these singular fibres. The
exceptional curves of this blow--up form the components of $C_+$. Because we
blow up points of multiplicity two on the fibres, the total transform of the
two singular fibres contains $2C_+$.  In the remaining two steps we have to
blow up smooth points on $A' + \overline{A}'$. Hence, we obtain $C_0 + 2C_+ \in
|2F|$. So, we can write $|-2K_S| = C_0 + 2C_- + |2F|$. Since we have
$\dim|\fund| = 1$, this is true for the generic real surface $S\in|\fund|$ by
Lemma \ref{sub} (b).
Let us denote by $\Phi = \Phi_{|-K|}$ the meromorphic map
$Z\dashrightarrow\bbfP^4$ defined by $|-K|$. If $\varphi:S\rightarrow\bbfP^2$
is the restriction  of $\Phi$ to a generic smooth real $S\in|\fund|$, then the
image of $\varphi$ is a conic in $\bbfP^2$. The general fibre of $\varphi$ is a
twistor fibre, hence a smooth rational curve intersecting $C$ transversally at
two points lying on $A$ and $\overline{A}$ respectively. 

Let $\tilde{Z}\rightarrow Z$ be a modification such that $\Phi$ becomes a
morphism $\tilde{\Phi}:\tilde{Z}\rightarrow\bbfP^4$. Because the smooth real
fundamental divisors $S$ sweep out a Zariski dense subset of $Z$, the image of
this set is also Zariski dense in $\tilde{\Phi}(\tilde{Z})\subset\bbfP^4$. As
the general fibre of $\Phi$, restricted to such surfaces $S$, is
one--dimensional, we obtain $\dim \tilde{\Phi}(\tilde{Z}) = 2$. Since 
$\tilde{Z}\rightarrow Z$ is a modification, there exists an
open Zariski dense subset $U\subset\tilde{\Phi}(\tilde{Z})$ such that the
fibres of $\tilde{\Phi}$ are irreducible curves. Moreover, we can choose $U$
such that the fibres of $\tilde{\Phi}$ over $U$ are isomorphic to
$\bbfP^1$,because this is true over a Zariski dense subset of
$\tilde{\Phi}(\tilde{Z})$.  Let $\tilde{\Phi}:\tilde{Z}_U\rightarrow U$ denote
the restriction of $\Phi$ over $U$. Then the preimage in $\tilde{Z}$ of the two
components of $C_-$ defines a pair of divisors $\Sigma$ and $\bar{\Sigma}$ in
$\tilde{Z}_U$ which are sections of $\tilde{Z}_U\rightarrow U$. Therefore,
${\cal E}:=\tilde{\Phi}_\ast\hol_{\tilde{Z}_U}(\Sigma + \bar{\Sigma})$ is a
vector bundle of rank three on $U$ and the canonical morphism
$\tilde{\Phi}^\ast{\cal E}\rightarrow \hol(\Sigma + \bar{\Sigma})$ is
surjective. This means that we obtain a morphism
$\tilde{Z}_U\rightarrow\bbfP({\cal E})$ which is compatible with the
projections to $U$. Restricted to each fibre this morphism is the Veronese
embedding of degree two $\bbfP^1\hookrightarrow\bbfP^2$. Hence, the image of
$\tilde{Z}_U$ in the quasi--projective variety $\bbfP({\cal E})$ is
three--dimensional. This implies $a(Z) = a(\tilde{Z}_U) = 3$, which completes
the proof.\qed
%
%
%
\section{Conclusions}
\label{vier}

In this section we collect the results of this paper to obtain a clear
picture of the situation considered.
By $Z$ we always denote a simply connected compact twistor space of positive
type over $4\PP$.

We call ${\mathcal N}:= \{C\subset Z\mid C$ irreducible curve, $C.(\fund)<0\}$
the set of negative curves. By definition, $\fdb$ is nef if and only if
${\mathcal N}\ne \emptyset$. The structure of ${\mathcal N}$ is described by
the following
\begin{thm}\label{ncurves}
  If ${\mathcal N}\ne \emptyset$ this set consists of a finite number of smooth
  rational curves. More precisely, only the following cases are possible:
  \begin{itemize}
  \item[(a)] ${\cal N}$ contains a real member $C_0$. Then: ${\cal N} =
    \{C_0\}$ and $C_0(\fund) = -2$, $\dim|\fund| = 2$ and $a(Z) = 3$.
  \item[(b)] ${\cal N}$ contains a non--real member $A$ with $A.(\fund) =
    -2$. Then, ${\cal N} = \{ A, \overline{A}\}$, $\dim|\fund| = 3, a(Z) = 3$
    and $Z$ is a LeBrun twistor space.
  \item[(c)] Each member $A\in {\cal N}$ fulfills $A.(\fund) = -1$. Then
    $|{\cal N}| \in \{2, 4, 6\}$, $\dim|\fund| = 1$ and $a(Z) =3$.
  \end{itemize}
\end{thm}
\pf We have only to collect the results of Section \ref{drei}.
\qed

We can compute the algebraic dimension in the following way:
\begin{thm}\label{main}
  $a(Z) = 3 \iff \fdb$ is not nef;\\
  $a(Z) = 2 \iff \fdb$ is nef and $\exists m\geq 1: h^1(\fb{m}) \ne 0$;\\
  $a(Z) = 1 \iff \forall m\geq 1: h^1(\fb{m}) = 0$.
\end{thm}
\pf
This results from Proposition \ref{tau} and Theorems \ref{nef} and
\ref{ncurves}.\qed 

We can characterize Moishezon twistor spaces as follows:
\begin{thm}
  The following conditions are equivalent:
  \begin{enumerate}
  \item $a(Z) = 3$;
  \item $\fdb$ is not nef;
  \item there exists a smooth rational curve $C\subset Z$ with $C.(\fund) < 0$.
  \end{enumerate}
\end{thm}
\pf
Apply Theorems \ref{ncurves} and \ref{main}.\qed

\begin{rem}
  Remembering that, by Poon's theorem, $\fdb$ is big if and only if $Z$ is
  Moishezon, we obtain from the preceding theorem: the line bundle $\fdb$ is
  never nef and big (under our special assumptions).
\end{rem}

LeBrun twistor spaces are characterized (see \cite{Ku}, \cite{Po2}) by
the property to contain a pencil of divisors of degree one. We can give (for
the case $n=4$) two further characterizations:
\begin{thm}\label{cb}
  The following properties are equivalent:
  \begin{enumerate}
  \item $Z$ contains a pencil of divisors of degree one;
  \item $\dim|\fund| = 3$;
  \item there exists a smooth rational curve $A\subset Z$ with $A \ne \bar{A}$
    and $A.(\fund) = -2$.
  \end{enumerate}
\end{thm}
\pf
The implications (i)$\Rightarrow$(ii) and (i)$\Rightarrow$(iii) follow from the
Kurke--Poon theorem. The reverse implications follow from Theorem
\ref{ncurves}.\qed
\begin{thm}
  $a(Z) \geq \dim|\fund|$.
\end{thm}
\pf
This follows directly from Proposition \ref{tau} and Theorems \ref{nef} and 
\ref{ncurves}. \qed

If $|\fund|$ is not a pencil, we obtain the following nice result:
\begin{thm}
  If $\dim|\fund|\geq 2$, then:\\
  $a(Z) = 2 \iff \fdb$ is nef $\iff |\fund|$ does not have base points.
\end{thm}
\pf
The first equivalence results from the previous theorem and Theorem \ref{main}.
If $|\fund|$ does not have base points, $\fdb$ is necessarily nef. If $\fdb$ is
nef and $\dim|\fund|\geq 2$ we have seen in Theorem \ref{nef} that $|\fund|$ is
base point free.\qed
\begin{cor}
  $|\fund|$ is base point free $\Rightarrow a(Z) = 2$.
\end{cor}
\pf
This is immediate from the previous theorem, because a pencil $|\fund|$ has
always base points. \qed
\begin{rem}
  The reverse implication is not true, which follows from the existence theorem
  in \cite{CK}. There, twistor spaces with $a(Z) = 2$ and
  $\dim|\fund| = 1$ over $4\PP$ were constructed.
\end{rem}
%
%

\end{document}